\documentclass{birkjour}
\usepackage{amsmath,amssymb,amsfonts,amsthm,mathrsfs}
\usepackage{mathtools}

\usepackage{graphicx,color}

\renewcommand{\theequation}{\arabic{section}.\arabic{equation}}
\theoremstyle{plain}
\numberwithin{equation}{section}

\newtheorem{theorem}{Theorem}[section]
\newtheorem{proposition}[theorem]{Proposition}
\newtheorem{lemma}[theorem]{Lemma}
\newtheorem{corollary}[theorem]{Corollary}
\newtheorem{remark}[theorem]{Remark}

\def\im{{\rm Im\, }}

\def\kr{{\rm Ker\, }}
\def\diag{{\rm diag\, }}

\def\rank{{\rm rank\, }}

\def\ind{{\rm ind\,}}

\def\BC{{\mathbb C}}
\def\BD{{\mathbb D}}

\def\BR{{\mathbb R}}
\def\BT{{\mathbb T}}

\definecolor{green}{rgb}{.0,.6,.0}

\newcommand{\bpr}{{\noindent\textbf{Proof.}\ \ }}
\newcommand{\epr}{{\hfill $\Box$}}

\begin{document}
\title {Wiener-Hopf indices of unimodular functions on the  imaginary  axis}

\author[A.E. Frazho]{A.E. Frazho}
\address{Department of Aeronautics and Astronautics,
Purdue University,
West Lafayette, IN 47907, USA.}
\email{frazho@purdue.edu}
\author[A.C.M. Ran]{A.C.M. Ran}
\address{Department of Mathematics, Faculty of Science, VU Amsterdam, De Boelelaan 1111, 1081 HV Amsterdam, The Netherlands and Research Focus: Pure and Applied Analytics, North-West~University, Potchefstroom, South Africa}
\email{a.c.m.ran@vu.nl}
\author[F. van Schagen]{F. van Schagen}
\address{Department of Mathematics, Faculty of Science, VU Amsterdam,
De Boelelaan 1111, 1081 HV Amsterdam}
\email{f.van.schagen@vu.nl}

\address[FvS]{Department of Mathematics, Faculty of Science, VU Amsterdam,
De Boelelaan 1111, 1081 HV Amsterdam}



\begin{abstract}
This paper is concerned with the Wiener-Hopf indices of unimodular rational matrix functions on the imaginary axis. These indices play a role in the Fredholm theory for Wiener-Hopf integral operators. Our main result gives formulas for the Wiener-Hopf indices in terms of the matrices appearing in realizations of the factors in a Douglas-Shapiro-Shields factorization of the unimodular function. Two approaches to this problem are presented: one direct approach using operator theoretic methods, and a second approach using the Cayley transform which allows to use results for an analogous problem regarding unimodular functions on the unit circle and corresponding Toeplitz operators.
\end{abstract}

\maketitle



{\date{}}



\setcounter{equation}{0}
\section{Introduction}

Wiener-Hopf factorization of matrix valued functions plays an important role in determining the Fredholm properties of several classes of operators, such as singular and Wiener-Hopf integral operators and Toeplitz operators, see e.g., \cite{BS, CGOT3,GF,GGKOT49,GGKOT63} and \cite{GKS}. To make this more explicit, let $R(s)$ be an $m\times m$ rational matrix valued function on the imaginary axis, which is
continuous and takes invertible values for $s$ on the imaginary axis.
A factorization
\[ \textstyle{
R(s)=W_-(s) \, \diag \left( (\frac{1-s}{1+s} )^{\kappa_j}\right)_{j=1}^m W_+(s), }
\]
where $W_-$ and its inverse are analytic on the closed left half of the complex plane, including infinity,
and $W_+$ and its inverse are analytic on the closed 
right hand half plane, including infinity, and $\kappa_j\in\mathbb{Z}$ for $j=1, \ldots, m$,
is called a (\emph{right}) \emph{Wiener-Hopf factorization} with respect to the imaginary  axis.
The integers $\kappa_j$ are uniquely determined by $R$ and they are called the
\emph{Wiener-Hopf indices} of $R$.

In order to present the main results of this paper we recall some
definitions and fix some notations.
Throughout $\mathcal{E} $ is a finite dimensional Hilbert space and
 $L_+^2(\mathcal{E})$ denotes the space of
square integrable functions from $[0,\infty)$ to $\mathcal{E}$.
The Laplace transform $\mathfrak{L}$ is the unitary operator mapping
$L_+^2(\mathcal{E})$ onto $H^2(\mathcal{E})$ defined by
\begin{equation}\label{fourier}
(\mathfrak{L}f)(s) =\int_0^\infty e^{-st} f(t) dt
\qquad (f\in L_+^2(\mathcal{E})).
\end{equation}
Here $H^2(\mathcal{E})$ is the Hardy space of $\mathcal{E} $-valued functions
that are  analytic in the right hand half complex plane
$ \BC_+ = \{s\in \mathbb{C}:\Re(s) >0\}$
and square integrable on the imaginary axis.
Furthermore $H^\infty(\mathcal{E},\mathcal{E})$ is the Hardy space
of all functions $\Theta $ whose values are operators on $ \mathcal{E} $ and that
are analytic and uniformly bounded in the open right half plane $ \BC_+ $, i.e.,
\[
\|\Theta\|_\infty = \sup\{\|\Theta(s)\|: \Re(s)>0\} <\infty.
\]
(If $\mathcal{E} = \mathbb{C}$, then
$H^\infty(\mathbb{C},\mathbb{C})$ is denoted by $H^\infty$.)
Similarly, $H^\infty(\mathcal{E})$ is the Hardy space
consisting of the set of all $\mathcal{E} $-valued functions
$\Psi(s)$ that are analytic in the open right half plane and such that
\[
\|\Psi\|_\infty = \sup\{\|\Psi(s)\|: \Re(s)>0\} <\infty.
\]

Let $R(s)$ be a rational function taking unitary values on the finite dimensional space $\mathcal{E}$ for values of $s$ on the imaginary line. In particular, $R(\infty)$ is also a unitary operator. Let $r(t)$ be given via
$R(i\omega) = R(\infty) + (\mathfrak{L} r)(i \omega)$.

Define $T_R$ to be the Wiener-Hopf  operator on $L_+^2(\mathcal{E})$ determined by $R$, that is,
\begin{equation}\label{TWR}
(T_R f)(t) = R(\infty) f(t)+ \int_0^\infty r(t-\tau)f(\tau)d\tau
\qquad (f\in L_+^2(\mathcal{E})).
\end{equation}
For such an operator the image is denoted by $\im T_R$ and
the kernel or null space by $ \kr T_R$.

Throughout  $\zeta(s)=\tfrac{1-s}{1+s}$ is the conformal mapping which maps
the open right half plane $\mathbb{C}_+ $
onto the open unit disc $\mathbb{D} = \{z\in \mathbb{C}: |z| <1\}$.
Let $ -\kappa_1 , -\kappa_2  , \ldots, -\kappa_p  $ with
$ \kappa_1 \geq \kappa_2 \geq \cdots \geq \kappa_p  $ be the negative
Wiener-Hopf indices of the function $ R $.
Then (see \cite{GGKOT49}  Theorem XIII.3.2) we have that
the dimension $\mathfrak{n}(T_R) $ of the null space
$\kr  T_R  $ is given by
\[
\mathfrak{n}(T_R) = \sum_{\kappa_j \geq 1 } \kappa_j .
\]
Consider the function $ R $ multiplied by $ \zeta(s)^k $, which we denote by $ \zeta^k R $.
Define for $ k = 1, 2, \ldots $ the numbers $ \mu_k $ by
\[
\mu_k = \mathfrak{n}(T_{\zeta^{k-1} R} ) - \mathfrak{n}(T_{\zeta^{k} R} ) .
\]
Then (see Section 2 below)
\[
\kappa_j = \# \{ k \mid \mu_k \geq j \}.
\]
Here $\# E $ denotes the number of elements of the set $E$.

Finding the Wiener-Hopf indices in terms of matrices in a realization of the function $R$
is a problem that has already
some history, see \cite{BGKOT12,BGK4OT21,BGKOT21,BGKROT200,GKS,GLeR}.
An analogous problem for a unimodular function on the unit circle was studied in \cite{GKRa}.
There, significant use was made of the
Douglas-Shapiro-Shields factorization of $R$, that is,
writing $R$ as $R=VW^*$, where $V$ and $W$ are bi-inner.
Note that in this case we say that the operator valued function
$\Theta$ is a {\it bi-inner function} if
$\Theta$ is analytic and uniformly bounded on $ \BD$, and
$\Theta( e^{i \omega})$ is a unitary operator on $\mathcal{E}$
for almost all $ \omega \in \BR $.
The Wiener-Hopf indices were given in \cite{GKRa} in terms of realizations
of $V$ and $W$, based on earlier work in \cite{FK}. These results were
extended to formulas for the Wiener-Hopf indices for any rational
matrix valued function in \cite{GKRb}.

Our aim in this paper is to obtain a result analoguous to the result
of \cite{GKRa} but for the imaginary axis replacing the unit circle and
with a different method, which is more operator theoretic. This method
is in parallel with our earlier results in \cite{FKRvS} for the case of
unimodular rational matrix functions on the unit circle and their
related Toeplitz operators.
Our approach leads to simple formulas for the Wiener-Hopf
indices of a rational matrix function that takes unitary values
on the imaginary axis.

Finally in this introduction we give a short description of the
various sections of this paper.
In Section 2, we introduce the functions and their realizations,
and present the main result in Theorem \ref{mainthm00}.
Section 3 is concerned with the Wiener-Hopf and Hankel operators corresponding to
bi-inner rational matrix functions.
Section 4 gives more detailed results on the unimodular function $R$,
its Wiener-Hopf operator $T_R$ and factorization.
In subsection 4.1 we specify the results for the case when
the bi-inner functions are scalar valued Blaschke products.
We derive the main results for this special case.
Section 5 is devoted to the proof of the main theorem.
In Section 6 we treat the Cayley transform and the connection it gives between
the realizations of rational matrix functions on the unit circle in
the complex plane and the realizations of rational matrix functions on
the imaginary axis in the complex plane.
In Section 7 these relations are used to prove the equivalence of the main result,
Theorem 2.2, of the paper with \cite[Theorem 2.2]{FKRvS}.
In Subsection 7.1 we present an example.


\section{The main result.}
To present our method to compute the Wiener-Hopf indices,
let us fix some notation.
Recall that $\{A \mbox{ on } \mathcal{X},B,C,D\}$ is a {\it realization of a function
$\Theta(s)$} if
\[
\Theta(s) = D +C(sI-A)^{-1}B.
\]
Here $A$ is an operator on $\mathcal{X}$ and $B$ maps $\mathcal{U}$ into $\mathcal{X}$,
while $C$ maps $\mathcal{X}$ into $\mathcal{Y}$ and
$D$ maps $\mathcal{U}$ into $\mathcal{Y}$.
Two state space realizations $\{A \mbox{ on }  \mathcal{X},B,C,D\}$
 and $\{A_1 \mbox{ on } \mathcal{X}_1 ,B_1,C_1,D_1\}$
are \emph{unitarily equivalent}  if $D = D_1$ and there
exists a unitary operator $U$ mapping $\mathcal{X}_1$ onto $\mathcal{X}$ such that
\[
A U = UA_1 \quad \mbox{and}\quad B = UB_1
\quad \mbox{and}\quad CU = C_1.
\]
We say that an operator $A$ on $\mathcal{X}$ is \emph{dissipative} if
$A^* +A \leq 0$.
An operator $A$ on a finite dimensional space $\mathcal{X}$ is
\emph{stable} if all the eigenvalues for $A$ are contained in
the open left hand plane $\mathbb{C}_-=\{s \in \mathbb{C}:\Re(s) <0\}$.

Throughout we will be dealing with finite dimensional realizations, that is,
realizations of the form $\{A \mbox{ on } \mathcal{X} ,B,C,D\}$ where
the state space $\mathcal{X}$ is finite dimensional.
In general we will call such a realization MIMO (multi input, multi output).
In the particular case  when $ \mathcal{U}$ and $\mathcal{Y} $ are one-dimensional
we say that the realization is SISO (single input, single output).
The realization
$\{A,B,C,D\}$ is \emph{stable,} if  $A$ is stable.
Next, we say that the realization
$\{A,B,C,D\}$ is stable and {\it dissipative } if the following holds:
\begin{enumerate}
  \item The operator $A$ is stable and $A+A^* +C^*C =0$;
  \item the operator $D$ is unitary;
  \item $B = - C^*D$.
\end{enumerate}
It is noted that the previous three conditions are equivalent to
\begin{itemize}
  \item[a)] The operator $A$ is stable and $A+A^* +BB^* =0$;
  \item[b)] the operator $D$ is unitary;
  \item[c)] $C = - D B^*$.
\end{itemize}

We are now ready to
present the following classical result, compare, e.g., \cite{AG}, \cite{BGKROT200}, Section 17.5.
The result is closely related to results on so called Livsic-Brodskii characteristic
operator functions, see Section 1.2 in \cite{BGKROT178}.

\begin{theorem}\label{thm-inner}  Let $\Theta$ be a rational function in $H^\infty(\mathcal{E},\mathcal{E})$.
Then $\Theta$ is bi-inner if and only if $\Theta$ admits
a stable dissipative  realization $\{A,B,C,D\}$. In this
case, all stable dissipative  realizations of $\Theta$ are
unitarily equivalent.
\end{theorem}

Let $R$ take unitary values on the imaginary axis.
Due to the Douglas-Shapiro-Shields factorization, when computing the 
Wiener-Hopf indices of $R$, without loss of generality, one can assume that
$R = VW^*$ where $V$ and $W$ are two rational
bi-inner functions in $H^\infty(\mathcal{E},\mathcal{E})$.
Here $ W^* $ denotes the function defined by
$ W^*(s) = \bigl( W(-\bar{s}) \bigr)^* $ ( $ s \in \BC $) and thus we have that
$ W^*(i \omega)  = \left( W(i \omega) \right)^* $ for all $\omega \in \BR$.
Let  $\{A_v \mbox{ on } \mathcal{X}_v, B_v,C_v,D_v \}$
and $\{A_w \mbox{ on } \mathcal{X}_w, B_w,C_w,  D_w\}$  be
two stable dissipative realizations of $V$ and $W$ respectively.
In particular,
\begin{align}
V(s) &= D_v +  C_v (s I -  A_v )^{-1} B_v \label{defV00} \\
W(s) &= D_w +   C_w (s I -  A_w )^{-1} B_w. \label{defW00}
\end{align}
Let $T_R$ be the Wiener-Hopf  operator on $L_+^2(\mathcal{E})$ determined by $R$ given by \eqref{TWR}.

Let $ -\kappa_1 , -\kappa_2  , \ldots, -\kappa_p  $ with
$ \kappa_1 \geq \kappa_2 \geq \cdots \geq \kappa_p  $ be the negative
Wiener-Hopf indices of the function $ R $
and $ R(s) = W_-(s) D(s) W_+(s) $ the Wiener-Hopf factorization of $ R $.
Then (see \cite{GGKOT49}  Theorem XIII.3.2) we have that
the dimension $\mathfrak{n}(T_R) $ of the null space of $ T_R $, is given by
\[
\mathfrak{n}(T_R) = \sum_{\kappa_j \geq 1 } \kappa_j .
\]
Recall that $\zeta(s) = \frac{1-s}{1+s}$.
Consider the function $ R $ multiplied by $ \zeta^k $, which we denote by $ \zeta^k R $.
Since $ \zeta^k(s) R(s) = W_-(s) \bigl( \zeta^k(s) D(s) \bigr) W_+(s)$,
the Wiener-Hopf indices of $\zeta^k R $ are each $k$ higher  than
the corresponding  index of  $R$.
Therefore
\[
\mathfrak{n}(T_{\zeta^k R} ) = \sum_{\kappa_j \geq k+1 }  (\kappa_j -k).
\]
Define for $ k = 1, 2, \ldots $, the numbers $ \mu_k $ by
\begin{equation}\label{defmuk10}
\mu_k = \mathfrak{n}(T_{\zeta^{k-1} R} ) - \mathfrak{n}(T_{\zeta^{k} R} ) =
\# \{ j : \kappa_j \geq k \}.
\end{equation}
Then (see \cite{GKvS1995} Proposition III.4.1)
\begin{equation}\label{defkappaj}
\kappa_j = \# \{ k : \mu_k \geq j \}.
\end{equation}

\medskip
The main result of the paper is the following theorem.
\medskip

\begin{theorem}\label{mainthm00}
Assume that $R = V W^*$ where
 $V$ and $W$ are two bi-inner rational functions in $H^\infty(\mathcal{E},\mathcal{E})$.
Let $\{A_v, B_v,C_v,D_v \}$
and $\{A_w, B_w,C_w,  D_w\}$  be
stable dissipative realizations of $V$ and $W$, respectively.
Let $\Omega$ be the unique solution of the Lyapunov  equation
\begin{equation}\label{lyapom00}
  A_v \Omega +\Omega A_w^* + B_v B_w^*=0.
\end{equation}
Let $C_\circ$ be the operator mapping $\mathcal{X}_w$ into $\mathcal{E}$ defined by
\begin{equation}  
C_\circ = D_v B_w^* + C_v \Omega.
\end{equation}
Finally,  let $Q$ be the unique solution to the Lyapunov  equation
\begin{equation}\label{defQ00}
 A_w  Q +Q  A_w^* + C_\circ^* C_\circ =0.
\end{equation}
Then the following holds:
\begin{enumerate}
  \item The operator $Q$ is a positive contraction.
  \item The multiplicity of $1$ as an eigenvalue of $ Q $  equals $\mathfrak{n}(T_R)$.
In other words, $\mathfrak{n}(T_R) = \mathfrak{n}(I-Q)$.
Moreover, for $\zeta(s) = \frac{1-s}{1+s}$, we have
\begin{equation}\label{alphathm00}
\mathfrak{n}(T_{\zeta^k R})  = \dim \left( \kr (I - \zeta(-A_w)^{k} Q \zeta(-A_w)^{*k} )\right).
\end{equation}
  \item For $k=1,2,\cdots$, consider the integers
\begin{equation}\label{defmuk00}
\mu_k = \mathfrak{n}(I - \zeta(-A_w)^{k-1} Q (\zeta(-A_w)^*)^{k-1} )
- \mathfrak{n} (I - \zeta(-A_w)^{k} Q \zeta(-A_w)^{*k} ) .
\end{equation}
Then the negative Wiener-Hopf indices $ -\kappa_1 , \ldots , -\kappa_p $
of $ T_R $ are given by
\begin{equation}\label{defmuk001}
\kappa_j = \# \{ k : \mu_k \geq j \}, \quad (j = 1 , \ldots, p = \mu_1 ).
\end{equation}
\end{enumerate}
\end{theorem}

Notice that, once Parts 1 and 2 are proven,  Part 3 follows from
the equations \eqref{defmuk10} and \eqref{defkappaj}.

The dual statement for the positive Wiener-Hopf indices is obtained by applying
the above theorem to the function $ R^*(s)=   R(-\overline{s})^*$;
see Corollary~\ref{maincol} below.


\setcounter{equation}{0}
\section{Unitary functions on the imaginary axis. }

Let $ \Theta \in H^\infty(\mathcal{E},\mathcal{E})$ be such that
$\Theta(s) - \Theta(\infty)= \left(\mathfrak{L}\theta\right)(s)$ and
$  \left(\mathfrak{L}\theta\right)(s) $ is
the  Laplace transform of  a rational function $\theta(t)$
whose values are linear operators on $\mathcal{E}$  and
such that $\int_0^\infty \|\theta(t)\|dt < \infty$.
(Here $\Theta(\infty)$ is the constant operator
function on $\mathcal{E}$ defined by
$\Theta(\infty) = \lim_{s\rightarrow \infty} \Theta(s)$.)
Then $T_\Theta$ is the \emph{Wiener-Hopf operator}  on $L_+^2(\mathcal{E})$
and $H_\Theta$ is the \emph{Hankel  operator} on $L_+^2(\mathcal{E})$
respectively defined by
\begin{equation}\label{deftoep}
(T_\Theta f)(t) = \Theta(\infty) f(t)+ \int_0^t \theta(t-\tau)f(\tau) d \tau
\end{equation}
\begin{equation} \label{defHank}
 (H_\Theta f)(t)  = \int_0^\infty \theta(t + \tau)f(\tau) d\tau
\end{equation}
where $f$ is a function in $L_+^2(\mathcal{E})$.
Throughout $\widetilde{\Theta}$ is the function in
$H^\infty(\mathcal{E},\mathcal{E})$ defined by
$\widetilde{\Theta}(s) = \Theta(\overline{s})^*$
for all $s$ in the open right half plane. In particular,
$\{A,B,C,D\}$ is a realization for $\Theta$ if and only if
$\{A^*,C^*,B^*,D^*\}$ is a realization for $\widetilde{\Theta}$.
Moreover, one readily obtains the following identity for Hankel operators:
$H_{\widetilde{\Theta}}=H_\Theta^*$.
Finally, it is noted that $ \theta(t) = C e^{A t} B$  for  $t\geq 0$.

As before, let $\Theta$ be a bi-inner function in
$H^\infty(\mathcal{E},\mathcal{E})$.
In this case, the corresponding Wiener-Hopf operator $T_\Theta$ is an isometry
on $L_+^2(\mathcal{E})$.
Indeed, $ T_{\Theta^\ast} T_\Theta = T_{\Theta^\ast \Theta } =  I $.
Because $\widetilde{\Theta}$ is also bi-inner,
$T_{\widetilde{\Theta}}$ is also an isometry on $L_+^2(\mathcal{E})$.
Let  $\mathfrak{H}(\Theta)$ and $\mathfrak{H}(\widetilde{\Theta})$
denote the orthogonal complements of the ranges of $T_\Theta$ and
$ T_{\widetilde{\Theta}} $, respectively, that is,
\begin{equation}\label{defHtheta}
\mathfrak{H}(\Theta)=L_+^2(\mathcal{E})\ominus T_\Theta L_+^2(\mathcal{E})
\quad \mbox{and}\quad
 \mathfrak{H}(\widetilde{\Theta}) = L_+^2(\mathcal{E})\ominus T_{\widetilde{\Theta}} L_+^2(\mathcal{E}).
\end{equation}
By consulting  Equation (24) in Section XII.2 of \cite{GGKOT49} and
with the fact that $\Theta$ is bi-inner, we see that
\begin{equation}\label{PH0}
P_{_{\mathfrak{H}( \Theta)}} = H_\Theta H_\Theta^* =
I - T_\Theta T_\Theta^*
\quad \mbox{and}\quad
P_{_{\mathfrak{H}(\widetilde{\Theta})}} = H_\Theta^* H_\Theta =
I - T_{\widetilde{\Theta}} T_{\widetilde{\Theta}}^*
\end{equation}
are orthogonal operators on $\mathfrak{H}({\Theta})$ and $\mathfrak{H}\widetilde{\Theta})$, respectively.
The second identity follows from the fact that $ \widetilde{\Theta} $ is also bi-iner.
In particular, $P_{_{\mathfrak{H}(\Theta)} } = H_\Theta H_\Theta^*$ and
the range of $H_\Theta$ equals $\mathfrak{H}(\Theta)$.
From the second equality in \eqref{PH0} we also have that
 $P_{_{\mathfrak{H}(\widetilde{\Theta})} } = H_\Theta^* H_\Theta$
and the range of the Hankel operator $H_\Theta^*$ equals $\mathfrak{H}(\widetilde{\Theta})$.
Therefore
\[
\im(H_\Theta) = \mathfrak{H}(\Theta) \quad \mbox{and}\quad
\kr(H_\Theta)^\perp = \mathfrak{H}(\widetilde{\Theta}).
\]

Using $P_{_{\mathfrak{H}(\widetilde{\Theta})} } = H_\Theta^* H_\Theta$, we see that
there exists a unitary operator $U$  mapping $\mathfrak{H}(\widetilde{\Theta})$ onto
$\im(H_\Theta) = \mathfrak{H}(\Theta)$ such that $U P_{_{\mathfrak{H}(\widetilde{\Theta})} } = H_\Theta$.
Therefore the Hankel operator $H_\Theta$ can be viewed as a unitary operator
mapping $\kr(H_\Theta)^\perp = \mathfrak{H}(\widetilde{\Theta})$ onto $\im (H_\Theta)=\mathfrak{H}(\Theta)$.
In particular, $\mathfrak{H}(\widetilde{\Theta})$ and  $\mathfrak{H}(\Theta)$
have the same dimension.

Let $\{A \mbox{ on } \mathcal{X},B,C,D\}$ be any stable dissipative
realization of a rational  bi-inner function $\Theta$ in $H^\infty(\mathcal{E},\mathcal{E})$.
Then its \emph{observability operator} $\Gamma $ mapping $\mathcal{X}$ into $L_+^2(\mathcal{E})$ and
\emph{controllability operator} $\Upsilon$ mapping $L_+^2(\mathcal{E})$ into $\mathcal{X}$ are
 defined by
\begin{align}\label{obsdef}
\Gamma x&= C e^{At}x \qquad \qquad \qquad (x\in \mathcal{X}, \ \ t\geq 0) \\
  \label{condef}
\Upsilon u &= \int_0^\infty e^{A t} B u(t) dt
                  \qquad (u\in L_+^2(\mathcal{E})).
\end{align}
Because $A^*+A + C^*C =0$ and $A$ is  stable,
the observability operator $\Gamma$ is an isometry.
Likewise, since $A^*+A + BB^* =0$ and $A$ is  stable,
the controllability  operator $\Upsilon $ is a co-isometry.
Using the fact that $\Theta(s) - \Theta(\infty)$
 is the Laplace transform of $ \theta(t):=C e^{At} B$,
 it follows that the Hankel operator $H_\Theta$ admits a
factorization of the form:
\begin{equation}\label{hankobs}
 H_\Theta = \Gamma \Upsilon.
\end{equation}
Recall that $H_\Theta$ can be viewed as a unitary
operator from $\mathfrak{H}(\widetilde{\Theta})$
onto $\mathfrak{H}(\Theta)$.
Since $\Gamma $ is an isometry and
$\Upsilon $ is a co-isometry, the equalities
$\mathfrak{H}(\Theta) = \im (\Gamma )$
and $\mathfrak{H}(\widetilde{\Theta}) = \im (\Upsilon^*)$ hold.
The equation $H_\Theta = \Gamma \Upsilon $ with \eqref{PH0},   readily
implies that
\begin{equation}\label{GammaH}
P_{_{\mathfrak{H}(\Theta)}} = \Gamma \Gamma^*
\quad \mbox{and}\quad
P_{_{\mathfrak{H}(\widetilde{\Theta})}} = \Upsilon^* \Upsilon.
\end{equation}


\setcounter{equation}{0}
\section{The function $R= VW^*$.  }
Let $V$ and $W$ be two rational bi-inner functions in
$H^\infty(\mathcal{E},\mathcal{E})$.
Let $R$ be the rigid function in $L^\infty(\mathcal{E},\mathcal{E})$ defined by
\begin{equation}\label{def-R}
R(i\omega) = V(i\omega) W(i\omega)^* \qquad (\mbox{for } -\infty < \omega < \infty).
\end{equation}
(A function $\Xi$ in $L^\infty(\mathcal{E},\mathcal{E})$ is \emph{rigid} if
$\Xi( i \omega )$ is almost everywhere a unitary operator on $\mathcal{\mathcal{E}}$.)
Since $V$, $W$ and $R$ are rational we extend the definition of $R$ to
all but a finite number of values of $ s \in \BC$ by
$R(s) = V(s) W(- \overline{s})^\ast  $.
Let $T_R$ be the Wiener-Hopf  operator on $L_+^2(\mathcal{E})$ determined by $R$,
see \eqref{TWR}.
Because $V$ and $W$ are bi-inner,  $H_V$ is a unitary operator from
$\mathfrak{H}(\widetilde{V})$ onto $\mathfrak{H}(V)$,
and $H_W$ is a unitary operator from
$\mathfrak{H}(\widetilde{W})$ onto $\mathfrak{H}(W)$.
Recall that $T_R = T_V T_W^* + H_V H_W^*$;
see Equation (24) in Section XII.2 of \cite{GGKOT49}.
This readily implies that
\begin{equation}\label{toephankY}
 T_R = T_V T_W^* + H_V Y H_W^*
\end{equation}
where $Y$ is the contraction mapping
$\mathfrak{H}(\widetilde{W})$ into
$\mathfrak{H}(\widetilde{V})$ defined by
\begin{equation}\label{defY}
Y = P_{_{\mathfrak{H}(\widetilde{V})}}|\mathfrak{H}(\widetilde{W})
: \mathfrak{H}(\widetilde{W}) \rightarrow \mathfrak{H}(\widetilde{V}).
\end{equation}
Since $V$ and $W$ are both bi-inner,
\[
L_+^2(\mathcal{E}) = \im (T_V) \oplus \im (H_V)
\quad \mbox{and}\quad
L_+^2(\mathcal{E}) = \im (T_W) \oplus \im (H_W).
\]
Using this with  $T_R = T_VT_W^* + H_VYH_W^*$, we see that
$T_R$ admits a "singular value type" decomposition of the form:
\begin{equation}\label{svdTR}
T_R = T_VT_W^* + H_VYH_W^*
= \begin{bmatrix}
    T_V & H_V \\
  \end{bmatrix}\begin{bmatrix}
                 I & 0 \\
                 0 & Y \\
               \end{bmatrix}\begin{bmatrix}
    T_W^* \\ H_W^* \\
  \end{bmatrix}.
\end{equation}
Here
\[
\begin{bmatrix}
    T_V & H_V \\
  \end{bmatrix}:\begin{bmatrix}
    L_+^2(\mathcal{E}) \\ \mathfrak{H}(\widetilde{V}) \\
  \end{bmatrix}\rightarrow L_+^2(\mathcal{E})
  \quad \mbox{and}\quad
  \begin{bmatrix}
    T_W^* \\ H_W^* \\
  \end{bmatrix}:L_+^2(\mathcal{E})
  \rightarrow \begin{bmatrix}
    L_+^2(\mathcal{E}) \\ \mathfrak{H}(\widetilde{W}) \\
  \end{bmatrix}
\]
are both unitary operators. Moreover, the
middle term
\begin{equation}\label{middle}
\begin{bmatrix}
                 I & 0 \\
                 0 & Y \\
               \end{bmatrix}  = \begin{bmatrix}
                 I & 0 \\
                 0 & P_{_{\mathfrak{H}(\widetilde{V})}}|\mathfrak{H}(\widetilde{W}) \\
               \end{bmatrix}:
\begin{bmatrix}
    L_+^2(\mathcal{E}) \\ \mathfrak{H}(\widetilde{W}) \\
  \end{bmatrix}
\rightarrow \begin{bmatrix}
    L_+^2(\mathcal{E}) \\ \mathfrak{H}(\widetilde{V}) \\
  \end{bmatrix}
\end{equation}
is a contraction.

Due to the decomposition of $T_R$ in \eqref{svdTR}, it follows that all the
properties such as invertibility and Fredholmness
of the operator $T_R$ are the same as those of the contraction $Y$.

It is noted that $x$ is in $\kr(Y)$ if and only if
$x$ is in $\mathfrak{H}(\widetilde{W})$ and $P_{_{\mathfrak{H}(\widetilde{V})}} x = 0$,
or equivalently, $x$ is in $\mathfrak{H}(\widetilde{W})$ and
$x$ is in $\im (T_{\widetilde{V}}) = \mathfrak{H}(\widetilde{V})^\perp$.
In other words,
\[
\kr(Y) = \im (T_{\widetilde{V}})\bigcap \mathfrak{H}(\widetilde{W})
\quad \mbox{and} \quad
\kr(Y^*) = \im (T_{\widetilde{W}})\bigcap \mathfrak{H}(\widetilde{V}),
\]
where the second equality follows from a similar argument.

Recall that an  operator $T$ mapping $\mathcal{X}$ into $\mathcal{Y}$
admits a Moore-Penrose inverse $T^{pinv}$  if the operator $T\vert\kr(T)^\perp$ mapping
$\kr(T)^\perp$ into the range of $T$ is invertible. In this case, the
Moore-Penrose inverse of $T$ is given by $T^{pinv}=\left(T\vert\kr(T)^\perp\right)^{-1}P_{\im (T)}$.
By consulting the form of $T_R$ in \eqref{svdTR},
we obtain the following result.

\begin{proposition}\label{prop-R}
Let $R=VW^*$ where $V$ and $W$ are both bi-inner functions in
$H^\infty(\mathcal{E},\mathcal{E})$. Moreover, let $Y$
be the contraction mapping $\mathfrak{H}(\widetilde{W})$ into
$\mathfrak{H}(\widetilde{V})$ defined by
$Y = P_{_{\mathfrak{H}(\widetilde{V})}}|\mathfrak{H}(\widetilde{W})$.
Then the following hold.
\begin{enumerate}
  \item The operator $T_R$ is invertible if and only if
  $Y$ is invertible. In this case,
  \begin{equation}\label{invTR}
    T_R^{-1} = T_W T_V^* + H_W Y^{-1}H_V^*.
  \end{equation}
  \item The subspaces  $\kr(T_R)$ and $\kr(Y)$ have the same dimension.
  In fact,
  \begin{equation}\label{kerTR}
    \kr(T_R) = H_W \kr(Y) \quad \mbox{and}\quad
    \kr(Y) = \im (T_{\widetilde{V}})\bigcap \mathfrak{H}(\widetilde{W}).
  \end{equation}
  \item The subspaces  $\kr(T_R^*)$ and $\kr(Y^*)$ have the same dimension.
  In particular,
  \begin{equation}\label{kerTR1}
    \kr(T_R^*) = H_V \kr(Y^*)
    \quad \mbox{and}\quad
    \kr(Y^*) = \im (T_{\widetilde{W}})\bigcap \mathfrak{H}(\widetilde{V}).
  \end{equation}
  \item The subspaces $\im (T_R)^\perp$ and $\im (Y)^\perp$
  have the same dimension. In fact,
  \begin{equation}\label{kerTR2}
   \im (T_R)^\perp =  \kr(T_R^*) = H_V \kr(Y^*) = H_V \,\im (Y)^\perp.
  \end{equation}
  \item The operator $T_R$ admits a Moore-Penrose restricted inverse if and only if
  $Y$ admits a Moore-Penrose restricted inverse. In this case,
  \begin{equation}\label{MPenTR}
    T_R^{pinv} = T_W T_V^* + H_W Y^{pinv}H_V^*.
  \end{equation}
\end{enumerate}
\end{proposition}


\subsection{The Blaschke product case}

In this section, to gain some insight into the general case,  we will study the
contraction $Y=P_{_{\mathfrak{H}(\varphi)}}|\mathfrak{H}(m)$
mapping $\mathfrak{H}(m)$ into $\mathfrak{H}(\varphi)$ when
$m$ and $\varphi$ are scalar Blaschke  products.
We say that a function $b(s)$ is a \emph{Blaschke product}  if
\begin{equation}\label{blaschke}
b(s) = \rho \prod_{k=1}^n \frac{s+\overline{\alpha}_k}{s- \alpha_k }
\qquad (\mbox{where } \Re(\alpha_k) < 0 \mbox{ for all } k).
\end{equation}
(Here $\rho$ is a complex number on the unit circle.)
Throughout $\widetilde{b}(s) = \overline{b(\overline{s})}$.
Finally, it is noted that $\widetilde{b}(-s) = \frac{1}{b(s)}$.

Moreover, $n = \deg(b)$ is the \emph{degree}  of the Blaschke product.
We will only consider Blaschke products of finite degree.
So if we say that $b(s)$ is a Blaschke product,
then we assume that  $b(s)$  is a function of the form \eqref{blaschke}
and the degree of $b$ is finite. It is well known that
$b(s)$ is a rational inner function in $H^\infty$ if and only if
$b$ is a Blaschke product (of finite degree).
Furthermore, the SISO (single input single output) function  $b(s)$
is a Blaschke product of degree $n$
if and only if $b(s)$ admits a stable dissipative realization
 $\{A \mbox{ on } \mathcal{X}, B ,C ,D\}$ where
$n$ is the dimension of the state space $\mathcal{X}$.
(See, eg., \cite{BGKROT200} Section 17.5.)
In this case, the poles of $b(s)$ are precisely the
eigenvalues of $A$.

Now assume that $\{A \mbox{ on } \mathcal{X}, B ,C ,D\}$  is a
stable dissipative realization for a Blaschke product $b(s)$
of degree $n$.
Let $\Gamma$ mapping $\mathcal{X}$ into $L_+^2$ be the
observability operator formed by the pair $\{C,A\}$.
Recall that $\Gamma$ is an isometry. Moreover, the
range of $\Gamma$ equals $\mathfrak{H}(b)$.
Because the dimension of the state space is $n$, it follows
that the dimension of $\mathfrak{H}(b)$
(denoted by $\dim(\mathfrak{H}(b))$) equals $n$.
Furthermore, the Laplace transform of the space $\mathfrak{H}(b)$ is given by
\begin{equation}\label{FHb}
\mathfrak{L}\big(\mathfrak{H}(b)\big) = \{\mathfrak{L}\big( \Gamma x\big) :x\in \mathcal{X}\}
= \{C(sI -   A)^{-1} x : x \in \mathcal{X}\}.
\end{equation}
Here $\mathfrak{L}\big(\mathfrak{H}(b)\big) $ stands for the set of
the Laplace transforms of the elements in $ \mathfrak{H}(b) $.
Because the pair $\{C,A\}$ is observable,
\begin{equation}\label{FHb1}
\mathfrak{L}\big(\mathfrak{H}(b)\big) = \left\{\frac{p(s)}{\det[sI- A]}: p(s) \mbox{ is a polynomial
of degree} < \dim({\mathcal{X}})\right\}.
\end{equation}
(To see this, if $f(s) \in\mathfrak{L}\big(\mathfrak{H}(b)\big)$,  then
$ f(s) = \frac{p(s)}{\det[sI- A]} $, where $ p(s)  $
is a polynomial of degree less than $\dim({\mathcal{X}}) $.
Moreover, the linear spaces in the right hand sides of \eqref{FHb} and \eqref{FHb1}
both have dimension $ n$. Therefore they are equal.)
In fact, $\det[sI- A] = \prod_{1}^n(s- \alpha_k)$.
So if $b(s)$ is the Blaschke product of degree $n$ given in
\eqref{blaschke}, then
\begin{equation}\label{FHb2}
\mathfrak{L} \big(\mathfrak{H}(b)\big) = \left\{\frac{p(s)}{\prod_{k=1}^n(s- \alpha_k)}: p(s) \mbox{ is a polynomial of degree} < n\right\}.
\end{equation}

Let $G(s) = G(\infty) + (\mathfrak{L} g)(s)$ be a function in $H^\infty$ where $g$ is in $L_+^1$.
Let $A$ be a stable operator  on $\mathcal{X}$.
Then
\begin{align*}
\left(T_G^* e^{A(\cdot)}\right)(t) &= \overline{G(\infty)}e^{At} + \int_t^\infty g(\tau-t)^*e^{A \tau} d \tau\\
&=\overline{G(\infty)}e^{At} + e^{At} \int_0^\infty g(v)^*e^{A v} d v\\
&=
\overline{G(\infty)}e^{At} + e^{At} \left[\int_0^\infty g(\tau)^*e^{-s \tau} d \tau\right]_{\mbox{evaluated at $-A$}} \\
&=
\overline{G(\infty)}e^{At} + e^{At}\big( \mathfrak{L}(g(\cdot)^*\big)(-A) \\
&= e^{At} \widetilde{G}(-A).
\end{align*}
Hence
\begin{equation}\label{eigen}
\big( T_G^* e^{A(\cdot)}\big)(t) =  e^{At}\widetilde{G}(\cdot)(-A).
\end{equation}

Let $A$ be a stable dissipative operator on $\mathcal{X}$. Let
$C$ be any operator from $\mathcal{X}$ onto $\mathcal{E}$
such that $A^* +A + C^*C =0$. Let $\Gamma$ be the observability operator
from $\mathcal{X}$ into $L_+^2(\mathcal{E})$ defined by the pair $\{C,A\}$.
Recall that $\Gamma$ is an isometry. Let $\psi(s)$ be any function in $H^\infty$.
Then
\[
\left(T_\psi^* \Gamma\right)(t) = C\Big( T_\psi^* e^{A(\cdot)}\Big)(t) = C e^{At}\widetilde{\psi}(-A)=
\big(\Gamma(t)\big) \widetilde{\psi}(-A).
\]
In other words,
\begin{equation}\label{funcpsi}
 T_\psi^* \Gamma = \Gamma \widetilde{\psi}(-A)
 \quad \mbox{and thus}\quad \Gamma^* T_\psi^* \Gamma = \widetilde{\psi}(-A).
\end{equation}
Therefore if $A$ is a stable dissipative operator and $\psi$ is a function
in $H^\infty$, then the evaluation $\psi(-A)$ of the function $\psi(s)$ at
$-A$ is given by
 \begin{equation}\label{funcpsi5}
  \psi(-A)  =  \Gamma^* T_{\widetilde{\psi}}^* \Gamma.
\end{equation}
Using this along with (again) the fact that $\Gamma$ is an isometry, we obtain
\[
\|\psi(-A)\| = \|\Gamma^* T_{\widetilde{\psi}}^* \Gamma \|
\leq \|T_\psi^*\|  =\|\psi\|_\infty.
\]
This readily implies that
\begin{equation}\label{funcA}
\|\psi(-A)\| \leq \|\psi\|_\infty \qquad (\mbox{$\psi\in H^\infty$ and $A$ is stable and dissipative}).
\end{equation}
In particular, if $\psi$ is an inner function, then
$\psi(-A)$ is a  contraction.
Finally, if $\psi$ and $\theta$ are two functions in $H^\infty$, it follows
that
\[
\big(\theta\psi\big)(-A) = \theta(-A)\psi(-A)
 \]
 Since $\psi$ and $\theta$ commute, $\theta(-A)\psi(-A) = \psi(-A)\theta(-A)$.

Recall that if $T$ is a contraction mapping $\mathcal{X}$ into $\mathcal{Z}$, then
 $D_T$ is the positive square root of $I-T^*T$, and $\mathfrak{D}_T$
 is the (closed) range of $D_T$. Finally, $\mathfrak{d}_T$ is the dimension of $\mathfrak{D}_T$.
 This sets the stage for the following result.

\begin{lemma}\label{lem-mphi} Let $m(z)$  and $\varphi(z)$ be  two finite  Blaschke products
in $H^\infty$.
Consider the contraction $Y=P_{_{\mathfrak{H}(\varphi)}}|\mathfrak{H}(m)$
mapping $\mathfrak{H}(m)$ into $\mathfrak{H}(\varphi)$.
Let $\{A,B,C,D\}$ be a   stable dissipative   realization for $m$.
Then the following holds.
\begin{enumerate}
  \item There exists a unitary operator $\Psi$ mapping the range of
  $Y$ onto $\mathfrak{D}_{_{\widetilde{\varphi}(-A)}}$ such that
   \begin{equation}\label{defectPsi}
     \Psi P_{_{\mathfrak{H}(\varphi)}}\Gamma_m = D_{_{\widetilde{\varphi}(-A)}}.
   \end{equation}
   (Here $\Gamma_m$ is the observability operator formed by the state space
   realization $\{A \mbox{ on } \mathcal{X},B,C,D\}$ for $m$.)
   In this case,
   \begin{equation}\label{defect00}
    \mathfrak{d}_{\widetilde{\varphi}(-A)} =  \mathfrak{d}_{ \varphi (-A)} =\min\left\{\deg(\varphi),\deg(m)\right\}  =
    \rank\big( P_{_{\mathfrak{H}(\varphi)}}\vert \mathfrak{H}(m)\big).
    \end{equation}
  \item In particular,  $\deg(\varphi) \leq \deg(m)$ if and only if  the range of
   the contraction $Y=P_{_{\mathfrak{H}(\varphi)}}|\mathfrak{H}(m)$ equals
   $\mathfrak{H}(\varphi)$. In this case,
   $\mathfrak{d}_{\varphi(-A)} =
   \mathfrak{d}_{\widetilde{\varphi}(-A)} =\deg(\varphi)$, and
   \begin{equation}\label{blaschke00}
   \dim(\kr(Y)) = \deg(m)-\deg(\varphi) \qquad (\mbox{when } \deg(\varphi) \leq \deg(m)).
   \end{equation}
  \item The operator $Y$ is one to one if and only if
  $\deg(m) \leq \deg(\varphi)$. In this case,
  $\mathfrak{d}_{\varphi(-A)} =
   \mathfrak{d}_{\widetilde{\varphi}(-A)} =\deg(m)$,
   \begin{equation}\label{dimmphi}
    \dim(\kr(Y^*)) = \dim(\im (Y)^\perp) = \deg(\varphi) - \deg(m).
\end{equation}
  \item If $\deg(\varphi) < \deg(m)$, then the Blaschke product
  \begin{equation}\label{aak}
   \varphi(s) = \frac{C(sI- A)^{-1}x}{C(sI- A)^{-1}\widetilde{\varphi}(-A)x}
    \qquad (\mbox{if } 0 \neq x \in \mathfrak{D}_{_{\widetilde{\varphi}(-A)}}^\perp).
  \end{equation}
\end{enumerate}
\end{lemma}

Recall that for a Blaschke product $\varphi$, we have
$\widetilde{\varphi}(-s) = \frac{1}{\varphi(s)}$. Hence
\begin{equation}\label{invb}
T_\varphi^* \Gamma_m = \Gamma_m \widetilde{\varphi}(-A)= \Gamma_m \varphi^{-1}(A).
\end{equation}

\noindent {\sc Proof.}
Because  $\{A \mbox{ on } \mathcal{X},B,C,D\}$ is a stable
dissipative  realization of $m$, we have $\deg(m) = \dim(\mathcal{X})$.
Recall that the Hankel operator
$H_m = \Gamma_m \Upsilon_m$ where $\Gamma_{m}$
maps $\mathcal{X}$ into $L_+^2$ is the observability
operator formed by $\{C,A\}$. Moreover, $\Upsilon_m $ mapping $L_+^2$
onto $\mathcal{X}$
is the controllability operator determined  by $\{A,B\}$. Furthermore,
$\Gamma_{m}$ is an isometry and $\Upsilon_m $ is a co-isometry.
Since the range of $H_m$ equals $\mathfrak{H}(m)$, it follows
that the subspace $\mathfrak{H}(m)$ equals the range of $\Gamma_{m}$.
Notice that $h$ is in $\mathfrak{H}(m)$ if and only if
$h = \Gamma_{m} x$ for some $x$ in $\mathcal{X}$.
In fact, this $x$ is uniquely determined by $h$ and given by
$x = \Gamma_{m}^* h$.
Using  $P_{_{\mathfrak{H}(\varphi)}}  = I - T_\varphi T_\varphi^*$
with $T_\varphi^* \Gamma_{m} = \Gamma_{m} \widetilde{\varphi}(-A)$
and $ T_\varphi $ is an isometry, we have
\begin{align}\label{long}
\|P_{_{\mathfrak{H}(\varphi)}} \Gamma_{m} x\|^2 &=
\|(I - T_\varphi T_\varphi^*)\Gamma_{m} x\|^2 =
\|\Gamma_{m} x\|^2 - \|T_\varphi T_\varphi^*\Gamma_{m} x\|^2\nonumber\\
&= \|x\|^2 - \| T_\varphi^*\Gamma_{m} x\|^2 =
\|x\|^2 - \|\Gamma_{m} \widetilde{\varphi}(-A) x\|^2\nonumber\\
&=\|x\|^2 - \|\widetilde{\varphi}(-A) x\|^2 =
 \langle x,(I - \widetilde{\varphi}(-A)^*\widetilde{\varphi}(-A))x \rangle   \nonumber\\
&= \|(I - \widetilde{\varphi}(-A)^*\widetilde{\varphi}(-A))^{\frac{1}{2}}x\|^2
=\|D_{_{\widetilde{\varphi}(-A)}}x\|^2.
\end{align}
Hence  there
exists a unitary operator $\Psi$ mapping the range of
  $Y$ onto $\mathfrak{D}_{_{\widetilde{\varphi}(-A)}}$ such that
   \begin{equation}\label{defectPsi54}
     \Psi P_{_{\mathfrak{H}(\varphi)}}\Gamma_{m} = D_{_{\widetilde{\varphi}(-A)}}.
 \end{equation}
This proves equation \eqref{defectPsi} in Part 1.

Now let us show  that $\mathfrak{d}_{_{\widetilde{\varphi}(-A)}} = \min\{\deg(\varphi),\deg(m)\}$.
To this end, first assume that $\deg(\varphi) \leq \deg(m)$.
Then we claim that $Y$ is onto $\mathfrak{H}(\varphi)$, and thus,
the rank of $Y$ equals  $\dim(\mathfrak{H}(\varphi)) = \deg(\varphi)$.
Assume that a vector $h\in \mathfrak{H}(\varphi)$ is orthogonal to the range of
$Y=P_{_{\mathfrak{H}(\varphi)}}|\mathfrak{H}(m)$.
Then it follows from
\eqref{kerTR1} with $ \widetilde{V} = \varphi $ and $ \widetilde{W}=m $ that
$h$ is also a vector in the range of $T_m$, that is,
$h\in \mathfrak{H}(\varphi)\cap T_m L_+^2 $.
By consulting \eqref{FHb1} or \eqref{FHb2}, we see that the Laplace transform
$\mathfrak{L}\big(\mathfrak{H}(\varphi)\big)$ of the subspace
$\mathfrak{H}(\varphi)$  consists of a set of rational functions,
with at most $\deg(\varphi)-1$ zeros. The Laplace transform of $h$ is given by
\[
\widehat{h}(s) =  \big(\mathfrak{L} h\big)(s) \in
\Big(\mathfrak{L}\big(\mathfrak{H}(\varphi)\big)\cap m H^2\Big).
\]
Since $\deg(\varphi) \leq \deg(m)$, and $m$ is a rational function  with
$\deg(m)$ zeros, the subspace
$\Big(\mathfrak{L}\big(\mathfrak{H}(\varphi)\big)\cap m H^2\Big) = \{0\}$.
Therefore $h=0$ and the operator $Y$ is onto, whenever
$\deg(\varphi)\leq\deg(m)$.
This with
$\Psi P_{_{\mathfrak{H}(\varphi)}}\Gamma_{m} = D_{_{\widetilde{\varphi}(-A)}}$, implies that
$\mathfrak{d}_{_{\widetilde{\varphi}(-A)}}= \deg(\varphi)$.
Replacing $\varphi$ with $\widetilde{\varphi}$ shows that
$\mathfrak{d}_{_{ \varphi (-A)}}= \deg(\varphi)$ when $\deg(\varphi) \leq \deg(m)$.

\medskip

Now assume that $\deg(m) \leq \deg(\varphi)$.
 Clearly,  $Y$ and $Y^*$ have the same rank.
Notice that $Y^*$ is the contraction determined by
\[Y^*  = P_{_{\mathfrak{H}(m)}}|\mathfrak{H}(\varphi):
\mathfrak{H}(\varphi) \rightarrow \mathfrak{H}(m).
\]
So $Y^*$ has the same form as $Y$, except $m$ and $\varphi$
interchange places. By our previous analysis
 $\rank(Y^*) = \deg(m)$ and $Y^*$ is onto $\mathfrak{H}(m)$.  So $Y$ is one to one.
Recall that  $\Psi P_{_{\mathfrak{H}(\varphi)}} \Gamma_{m} = D_{_{\widetilde{\varphi}(-A)}}$.
Because  $Y$ is one to one, $D_{_{\widetilde{\varphi}(-A)}}$ must also be one to one.
Since $D_{_{\widetilde{\varphi}(-A)}}$ is one to one and $\dim(\mathcal{X}) = \deg(m)$,
we see that $\mathfrak{d}_{_{\widetilde{\varphi}(-A)}} = \deg(m)$.
This completes the proof of Part 1.

\medskip

To prove Part 2, we showed that
if $ \deg (\varphi) \leq \deg(m)$  then $ Y$ is onto $ \mathfrak{H} (\varphi ) $.
Moreover, $ \dim(\kr(Y)) + \dim(\im(Y)) = \dim\left( \mathfrak{H}(m) \right)$,
which proves \eqref{blaschke00}.
On the other hand if $\deg(\varphi) > \deg(m)$, then $\rank Y = \rank Y^* =\deg{m} $
and hence $Y$ is not onto $ \mathfrak{H} (\varphi ) $.

Part 3 is proven in the same way by replacing $Y $ by $Y^* $.

To establish Part 4, assume that  $\deg(\varphi) < \deg(m)$.
Then there exists a nonzero $x$ such that $D_{_{\widetilde{\varphi}(-A)}}x =0$, or
equivalently, $x = \widetilde{\varphi}(-A)^* \widetilde{\varphi}(-A)x$.
Using
\[
\Psi P_{_{\mathfrak{H}(\varphi)}}\Gamma_{_{m}}x = D_{_{\widetilde{\varphi}(-A)}}x=0,
\]
we have
$P_{_{\mathfrak{H}(\varphi)}} \Gamma_{m}x =0$.
By employing
$P_{_{\mathfrak{H}(\varphi)}} = I-T_\varphi T_\varphi^*$,  we obtain
\[
0 = P_{_{\mathfrak{H}(\varphi)}} \Gamma_m x
= (I-T_\varphi T_\varphi^*) \Gamma_m x =
\Gamma_m  x - T_\varphi  \Gamma_m  \widetilde{\varphi}(-A) x.
\]
In other words,
$\Gamma_m x = T_\varphi  \Gamma_m  \widetilde{\varphi}(-A) x$.
By taking the Laplace transform, we arrive at
\begin{equation}\label{scalar}
C(s I  -   A)^{-1}x = \varphi(s)C(sI- A)^{-1}\widetilde{\varphi}(-A) x.
\end{equation}
Since $\Gamma_{m} $ is one to one and $\widetilde{\varphi}(-A)x$ is nonzero,
$C(sI  -  A)^{-1}\widetilde{\varphi}(-A) x$ is a nonzero  function in  $H^2$.
Hence the rational function  $C(sI  -   A)^{-1}\widetilde{\varphi}(-A)x$
is nonzero.
Dividing \eqref{scalar} by $C(s I-  A)^{-1}\widetilde{\varphi}(-A) x$,
 yields the formula that we have been looking for, that is,
\[
\varphi(s) = \frac{C(s I-  A)^{-1} x}{C(sI-  A)^{-1}\widetilde{\varphi}(-A) x}.
\]
This completes the proof.
\epr


\begin{remark} Let $R$ be the rational rigid function in $L^\infty$ defined  by
$R(i\omega) = \varphi(i\omega) \overline{m(i\omega)}$,
where $\varphi$ and $m$ are two Blaschke products.
Let $T_R$ be the Wiener-Hopf operator on $L_+^2$ determined by $R$.
By consulting Proposition \ref{prop-R} and Lemma \ref{lem-mphi}, we readily obtain
the following.
\begin{enumerate}
\item The operator $T_R$ is invertible if and only if  $\deg(\varphi) = \deg(m)$,
      that is, the numbers of zeros and poles in $\BC_+$ of $R(s)$ are equal.
\item The kernel of  $T_R$ is nonzero  if and only if
      $\deg(\varphi) < \deg(m)$. In this case, the operator
      $T_R$ is onto $L_+^2$ and $\dim(\kr(T_R)) = \deg(m) - \deg(\varphi)$.
      The number of poles in $ \BC_+ $ of $R(s)$ is higher than
      the number of zeros.
\item The subspace $\im(T_R)^\perp$ is nonzero  if and only if
       $\deg(m) < \deg(\varphi)$. In this case, the range of the operator
       $T_R$ is closed, $\kr(T_R) = \{0\}$ and
       $\dim(\im(T_R)^\perp) = \deg(\varphi) - \deg(m)$.
\end{enumerate}
\end{remark}


 \begin{corollary}\label{cor-phi}
 Let $A$ be a stable dissipative operator  on a finite dimensional space $\mathcal{X}$
 such that $A^*+A$ has rank one,
 and $\varphi$ a rational Blaschke product in $H^\infty$.
 Then
 \begin{equation}\label{min00}
\mathfrak{d}_{_{\varphi(-A)}} = \mathfrak{d}_{_{\varphi(-A^*)}}
= \min\{\deg(\varphi),\dim(\mathcal{X})\}.
\end{equation}
 \end{corollary}

 \noindent {\sc Proof.}
Let $C$ be any operator mapping $\mathcal{X}$ into $\mathbb{C}$ such that
$A^*+A = -C^*C$. Set $B = -C^*$ mapping $\mathbb{C}$ into $\mathcal{X}$ and
   $D =1$. Then $\{A,B,C,D\}$ is a dissipative realization for
a Blaschke product $m(z)$ with degree $\dim(\mathcal{X})$; see
Theorem \ref{thm-inner}.
 Applying Lemma \ref{lem-mphi} yields \eqref{min00}.
  \epr

  \begin{proposition}\label{prop-clt00}
  Let $A$ be a stable dissipative operator on $\mathcal{X}$ such that the
  rank of $A^*+A$ equals $1$.
  Let $\varphi$ be a Blaschke product where
  $\deg(\varphi) < \dim(\mathcal{X})$. Then the
  inner function $\varphi$ is given by
  \begin{equation}\label{clt}
\varphi(s) = \frac{C(sI-  A)^{-1}\varphi(-A^*)x}{C(sI-  A)^{-1}x}
\qquad (\mbox{where } x=\varphi(-A^*)^*\varphi(-A^*)x \mbox{ and } x \neq 0).
\end{equation}
Here $C$ is any operator mapping $\mathcal{X}$ into $\mathbb{C}$ such that $A^*+A = -C^*C$.
If $\theta$ is any function in $H^\infty$ such that
$\theta(-A) = \varphi(-A)$ and $\|\theta\|_\infty \leq 1$, then $\theta(s) = \varphi(s)$.
  \end{proposition}

\noindent
{\sc Proof.} Formula \eqref{clt} will be derived from formula \eqref{aak}.
To this end, recall that if $Z$ is a contraction, then
$Z$ is a unitary operator from $\mathfrak{D}_Z^\perp$ onto
$\mathfrak{D}_{Z^*}^\perp$. Hence $\widetilde{\varphi}(-A)$ is a unitary operator
from $\mathfrak{D}_{_{\widetilde{\varphi}(-A)}}^\perp$ onto
$\mathfrak{D}_{_{\varphi(-A^*)}}^\perp$. In particular,
\[
\widetilde{\varphi}(-A): \mathfrak{D}_{_{\widetilde{\varphi}(-A)}}^\perp
\rightarrow \mathfrak{D}_{_{\varphi(-A^*)}}^\perp
\]
is unitary. So for $x\neq 0$ in  $\mathfrak{D}_{_{\widetilde{\varphi}(-A)}}^\perp$, we see that
$ y =  \widetilde{\varphi}(-A)x$ is in $\mathfrak{D}_{_{\varphi(-A^*)}}^\perp$ and
$x= \widetilde{\varphi}(-A)^*y$. By interchanging the roles of $x$ and $y$, we see
that \eqref{aak} yields  \eqref{clt} and vice-versa.


Let us show that if  $\theta$ is a function
$H^\infty$ such that $\theta(-A^*) = \varphi(-A^*)$ and  the $H^\infty$ norm   $\|\theta\|_\infty \leq 1$,
then $\theta(s) = \varphi(s)$.

Assume that $\theta$ is a function in
$H^\infty$ such that $\widetilde{\theta}(-A) = \widetilde{\varphi}(-A)$
 and $\|\theta\|_\infty \leq 1$.
By taking the adjoint, we see that $\theta(-A^*) = \varphi(-A^*)$.
In this case,
\[
T_\theta^* \Gamma = \Gamma \widetilde{\theta}(-A) =
 \Gamma \widetilde{\varphi}(-A).
 \]
Since $T_\varphi^* \Gamma = \Gamma \widetilde{\varphi}(-A)$, we also have
$T_\varphi^* \Gamma = T_\theta^* \Gamma$, or equivalently, \allowbreak
$ \Gamma^* T_\varphi = \Gamma^*T_\theta$. Multiplying by
$\Gamma$ on both sides, we obtain with $x$ as above
\[
\Gamma\Gamma^*T_\theta \Gamma x =  \Gamma\Gamma^*T_\varphi \Gamma x
= \Gamma   \varphi(-A^*) x.
\]
Because $\Gamma\Gamma^*$ is an orthogonal projection,
$\| \varphi(-A^*) x\|^2 = \|x\|^2$ and $T_\theta$ is a contraction,
we see that
\[
\|x\| \geq \|T_\theta \Gamma x\| \geq
\|\Gamma\Gamma^* T_\theta \Gamma x\|=\|\Gamma\varphi(-A^*) x\| = \|\varphi(-A^*) x \| =\|x\|.
\]
Therefore we have equality, and thus,
\[
T_\theta \Gamma x = \Gamma\Gamma^* T_\theta \Gamma x = \Gamma\varphi(-A^*) x
\qquad (\mbox{when } 0\neq x \in \mathfrak{D}_{\varphi(-A^*)}^\perp).
\]
By taking the Laplace transform of both sides, and using \eqref{clt}, we obtain
\[
\theta(s) = \frac{C(s I-  A)^{-1} \varphi(-A^*) x}{C(sI-  A)^{-1}x} = \varphi(s).
\]
Therefore $\theta(s) = \varphi(s)$. In other words, if
$ \theta $ is a function in $H^\infty$ such that
$\widetilde{\theta}(-A) = \widetilde{\varphi}(-A)$ and $\|\theta\|_\infty \leq 1$,
then $\theta(s) = \varphi(s)$. Replacing $\widetilde{\varphi}$ by $\varphi$
shows that if $ \theta $ is a function in $H^\infty$ such that
$ \theta (-A) =  \varphi (-A)$ and $\|\theta\|_\infty \leq 1$,
then $\widetilde{\theta}(s) = \widetilde{\varphi}(s)$, and thus,  $\theta(s) = \varphi(s)$.
\epr

\bigskip

Finally, it is noted that this result is an application    of
the Sz.-Nagy-Foias commutant lifting theorem; see Corollary 2.7
page 142 of \cite{FFGK} and is also deeply connected to
some of the results in  \cite{AAK}. See also \cite{FF}.

 \begin{proposition}\label{prop-phi}
 Let $A$ be a stable dissipative operator on a finite dimensional space
 $\mathcal{X}$ and $\varphi$ a Blaschke product. If
 $\deg(\varphi) \geq \dim(\mathcal{X})$, then
 the defect index for $\varphi(-A)$ equals $\dim(\mathcal{X})$,
 or equivalently,
 $\varphi(-A)^*\varphi(-A)$ does not have $1$ as an eigenvalue.
 \end{proposition}

\noindent {\sc Proof.}
 Let $C$ be any operator from $\mathcal{X}$ onto $\mathcal{E}$ such
 that $A^*+A = -C^*C$. Then the observability operator
 $\Gamma $ formed by $\{C,A\}$ is an isometry
 from $\mathcal{X}$ into $L_+^2(\mathcal{E})$.
 Recall that if $u(z)$ is any function in $H^\infty$, then
 $\|u(-A)\| \leq \|u\|_\infty$; see \eqref{funcA}.
 Because  $\varphi$ is an inner function, $\widetilde{\varphi}(-A)$ is a contraction.
  For $x$ in $\mathcal{X}$, we have
 \begin{align*}
&\  \|\big(I - \widetilde{\varphi}(-A)^*\widetilde{\varphi}(-A)\big)^{\frac{1}{2}}x\|^2
 = \langle\big(I - \widetilde{\varphi}(-A)^*\widetilde{\varphi}(-A)\big)x,x\rangle 
  \\
 &=\|x\|^2 - \|\widetilde{\varphi}(-A) x\|^2= \|\Gamma x\|^2 - \|\widetilde{\varphi}(-A) x\|^2\\
 & =
 \|\Gamma x\|^2 - \|T_{\varphi I} T_{\varphi I}^* \Gamma x\|^2= \|(I -  T_{\varphi I}T_{\varphi I}^*) \Gamma x\|^2=
 \|P_{_{\mathfrak{H}(\varphi I)}}\Gamma x\|^2.
 \end{align*}
 Hence there exists a unitary operator
 $\Psi$ from the
 range of $P_{_{\mathfrak{H}(\varphi I)}}\Gamma $
 onto the range of the defect operator
 $\big(I - \widetilde{\varphi}(-A)^*\widetilde{\varphi}(-A)\big)^{\frac{1}{2}}$
 such that
 \begin{equation}\label{defect}
 \Psi P_{_{\mathfrak{H}(\varphi I)}}\Gamma  =
 \big(I - \widetilde{\varphi}(-A)^*\widetilde{\varphi}(-A)\big)^{\frac{1}{2}}.
 \end{equation}
 In particular, the defect index $\mathfrak{d}_{_{\widetilde{\varphi}(A)}}$
 equals the rank of $P_{_{\mathfrak{H}(\varphi I)}}\Gamma $.

 Notice that $x$ is in the kernel of $P_{_{\mathfrak{H}(\varphi I)}}\Gamma $
 if and only if $x$ is a vector with  eigenvalue $1$ for
 $\widetilde{\varphi}(-A)^*\widetilde{\varphi}(-A)$. In this case,
 \[
 0 = P_{_{\mathfrak{H}(\varphi I)}}\Gamma  x =
 \Gamma  x - T_{\varphi I}T_{\varphi I}^* \Gamma  x.
 \]
 Hence
 \[
 \Gamma x =  T_{\varphi I} \Gamma  \widetilde{\varphi}(-A)x
 \qquad (\mbox{for } x \in \mathfrak{D}_{\widetilde{\varphi}(-A)}^\perp).
 \]
 By taking the Laplace transform
 \begin{equation}
 C(sI-  A)^{-1}x = \varphi(s) C(sI- A)^{-1}\widetilde{\varphi}(-A)x
 \qquad (\mbox{for } x \in \mathfrak{D}_{\widetilde{\varphi}(-A)}^\perp).
 \end{equation}
Because $\varphi(s)$ has all its zeroes in the open right half plane,
 and the eigenvalues of $A$ are in the open left half plane,
the numerator of $C(sI-  A)^{-1}x$ has at least
 $\deg(\varphi)$ zeros. However, the numerator of $C(sI-  A)^{-1}x$
 is a polynomial of degree at most $\dim(\mathcal{X}) -1 < \deg(\varphi)$.
 Therefore $x=0$, and
 $\widetilde{\varphi}(-A)^*\widetilde{\varphi}(-A)$
has no eigenvalue on the unit circle. The same argument
applies by replacing $\widetilde{\varphi}(s)$ by
$\varphi(s)$.  This completes the proof.
\epr

 \paragraph{Winding numbers in the scalar case}
Let $R$ be the rigid function in $L^\infty$ defined by
$R(i\omega) = \varphi(i\omega) \overline{m(i\omega)} =
 (\varphi m^* ) (i\omega)$,
where both $\varphi$ and $m$ are Blaschke products.
The {\it winding number} for $R$ is given by $\deg(\varphi) - \deg(m)$.
As before, let $Y$ be the contraction from
$\mathfrak{H}(m)$ into $\mathfrak{H}(\varphi)$ given by
$Y = P_{_{\mathfrak{H}(\varphi)}}\vert \mathfrak{H}(m)$.
By consulting Lemma \ref{lem-mphi}, we see that
\begin{enumerate}
  \item The winding number for $R = \varphi m^*$ equals
  $\dim\big(\im(Y)^\perp\big)$ when\\ $\deg(m) \leq \deg(\varphi)$.
  \item The winding number for $R = \varphi m^*$ equals
  $-\dim\big(\kr(Y)\big)$ when\\ $\deg(\varphi) \leq \deg(m)$.
  \item  Hence   the winding number for  $R = \varphi m^*$ equals $-\ind(Y)$
 where $\ind(Y)$ denotes  the Fredholm index
  ($\dim(\kr(Y)) - \dim(\kr(Y^*))$)   for $Y$.
\end{enumerate}

Recall that the Wiener-Hopf  operator $T_R$ admits a  decomposition of the form:
\[
T_R = T_\varphi T_m^* + H_\varphi Y H_m^*;
\]
see \eqref{svdTR}.
Here $Y$ is the finite dimensional contraction from $\mathfrak{H}(\widetilde{m})$ into
$\mathfrak{H}(\widetilde{\varphi})$ determined  by
$Y = P_{_{\mathfrak{H}(\widetilde{\varphi})}}\vert \mathfrak{H}(\widetilde{\varphi})$.
Hence
$T_R$ is Fredholm. Moreover, $T_R$ and $Y$ have the
same Fredholm index.   Hence
\begin{align}\label{frecind}
\ind(T_R) &= \dim(\kr (T_R)) - \dim(\kr(T_R^*))
 = \ind(Y) =\deg(m) - \deg(\varphi).
\end{align}
 The Fredholm index of $T_R$ equals minus the
winding number of $R = \varphi m^*$.

We say that a polynomial $p(s)$ is \emph{stable}  if all the roots of
$p(s)$ are contained in the open left hand plane $\BC_-$.
Let $p^\sharp(s)$ be the polynomial defined by
\[
p^\sharp(s) = \overline{p(-\overline{s})}.
\]
It is noted that if $p(s)$ is a stable polynomial of
degree $n$, then $\varphi(s)= \frac{p^\sharp(s)}{p(s)}$ is a Blaschke
product of order $n$. In particular, if $A$ is a stable dissipative
operator and $A^* + A$ has rank one, then
\[
\varphi(-A) = p^\sharp(s)(-A) \big(p(-A)\big)^{-1}
\]
is a strict contraction; see Corollary \ref{cor-phi}. In other words,
if $p(s)$ is a stable polynomial of degree $n$, then
\[
p(-A)^* p(-A) - \big(p^\sharp(-A)\big)^*p^\sharp(s)(-A)
\]
is a strictly  positive.

The following is a continuous time version of
the Schur-Cohen stability test.
\begin{remark}\label{rem-SC} 
Assume that $p(s)$ is a polynomial of degree $n$ and let
\[
R(s) = \frac{p^\sharp(s)}{p(s)}.
\]
Then the following statements are equivalent.
\begin{enumerate}
  \item The polynomial $p(s)$ is stable.
  \item The winding number for $R(s)$ equals $n$.
  \item The Fredholm index for $T_R$ equals $-n$.
  \item If $A$ is a stable dissipative operator on
  a $n$ dimensional space and the rank of $A+A^*$
  equals $1$, then
  \[
    p(-A)^* p(-A) - \big(p^\sharp(-A)\big)^*p^\sharp(-A)
  \]
is a strictly  positive.
\end{enumerate}
\end{remark}

\bpr
First observe that if
$q$ is a stable polynomial of degree strictly less than $n$
and $\theta(s) = \frac{q^\sharp(s)}{q(s)}$, then the norm of
$\theta(-A)$ equals $1$, or equivalently,
\[
    q(-A)^* q(-A) - \big(q^\sharp(-A)\big)^*q^\sharp(-A)
\]
is positive and singular; see Corollary \ref{cor-phi}.

Now let $p(s)$ be any polynomial of degree $n$.
Then $p(s)$ admits a factorization of the form:
$p(s) = q(s) u(s)$ where $q(s)$ is stable and $u(s)$ is
unstable. Then
\begin{align*}
&p(-A)^* p(-A) - \big(p^\sharp(-A)\big)^*p^\sharp(-A) =\\
&
u(-A)^* q(-A)^* q(-A)u(-A)  - \big(q^\sharp(-A)\big)^*\big(u^\sharp(-A)\big)^*q^\sharp(-A)u^\sharp(-A)
\end{align*}

Next we prove that item (1) implies item (4).
Since $ R $ is a Blaschke product we know that $ R(-A) $ is a contraction.
(See text after formula (4.19)).
Moreover equality (4.28) shows that the defect index of $R(-A) $ is equal to $n$.
Thus we have that $ R^*(-A) R(-A) < I $ or $ I - R^*(-A) R(-A) > 0 $.
Multiplying the last inequality by $  p^*(-A) p(-A) $ gives that
\[
    p(-A)^* p(-A) - \big(p^\sharp(-A)\big)^*p^\sharp(-A)
  \]
is a strictly  positive.

To prove that item (4) implies item (1) assume that $p$ is not stable and
$  p(-A)^* p(-A) - \big(p^\sharp(-A)\big)^*p^\sharp(-A) >0 $.
We will show that this leads to an obvious contradiction.
Write $ p = p_1 p_2^\sharp $ with $ p_1 $ and $ p_2 $ both stable polynomials.
The assumption gives that $ p_2 $ is not a constant.
Note that $ p(-A) $ must be invertible.
It follows that the zeros of $p$ and the eigenvalues of $-A$ are disjunct sets.
Write
\[
B_i(s) = p_i^\sharp (s) p_i(s)^{-1} \quad (i=1,2),  \mbox{ and } R(s)= B_1(s) B_2(s)^{-1}.
\]
Since $ B_i$  is a Blaschke product, we have that $ B_i(-A) $ is a contraction.
Also our assumption gives that $ B_1(-A)^* B_1(-A) < B_2(-A)^* B_2(-A) $.
As the degree of $ B_1 $ is strictly less than $n$, also
the defect index, which is the minimum of $n$ and the degree of $ B_1$,
is strictly less than $n$.
Therefore there exists a vector $ v$ such that
\begin{equation}\label{defectB1}
 v^* B_1(-A)^* B_1(-A) v = v^* v .
\end{equation}
On the other hand
\begin{equation}\label{B2largerB1}
 v^* B_1(-A)^* B_1(-A) v < v^* B_2(-A)^* B_2(-A) v \leq v^* v .
\end{equation}
The two equalities \eqref{defectB1} and \eqref{B2largerB1} contradict and therefore $p_2 $ is a constant.
We conclude that $p$ is stable.
\epr

A discrete time version of part of  Remark \ref{rem-SC}
with a different type proof is presented in \cite{PY}. 

Finally, for an example of a stable dissipative  matrix $A$ 
 such that $A+A^*$ has rank one, consider
the lower triangular Toeplitz matrix
$A$ on $\mathbb{C}^n$ given by
\begin{align*}
A &= \begin{bmatrix*}[r] 
-1       & 0            & 0             & 0            & \cdots & 0  \\
       2        & -1           & 0             & 0            & \cdots & 0  \\
       -2       &  2           & -1            & 0            & \cdots & 0  \\
       2        & - 2          & 2             & -1           & \cdots & 0  \\
       \vdots   & \vdots       & \vdots        & \vdots       & \vdots & \vdots  \\
       (-1)^{^n} 2  & (-1)^{^{n-1}}2  &  (-1)^{^{n-2}}2  & (-1)^{^{n-3}}2  & \cdots & -1  \\
     \end{bmatrix*}  \\
C&= \sqrt{2}\begin{bmatrix}
              1 & -1 & 1 & -1 & \cdots  & (-1)^{n-1} \\
            \end{bmatrix}.
\end{align*}
Notice that $A$ is the lower triangular Toeplitz
matrix formed by \[
\{-1,2,-2,2,-2,\cdots, 2(-1)^n\}.
\]
Then $A$ is a stable dissipative matrix and the rank of $A+A^*$ equals one.
In fact, $A+A^*+C^*C =0$.  This matrix is
motivated by equation \eqref{need}
below.

\setcounter{equation}{0}
\section{Point evaluation and MIMO systems}

Let us return to the MIMO case
and provide a proof of one of our main results, Theorem \ref{mainthm00}.

As before, assume that $V$ and $W$ are two bi-inner rational functions
in $H^\infty(\mathcal{E},\mathcal{E})$.
Let $\{A_v \mbox{ on } \mathcal{X}_v, B_v,C_v,D_v\}$
and $\{A_w \mbox{ on } \mathcal{X}_w, B_w,C_w,D_w\}$
be stable dissipative  realizations of $V$ and $W$, respectively.
Recall that  $R=VW^*$ and the Wiener-Hopf operator  $T_R$ on $L_+^2(\mathcal{E})$ admits
a decomposition of the form:
\begin{equation}\label{TRdef8}
T_R = T_V T_{_W}^* + H_V Y H_{W}^*.
\end{equation}
Here $Y$ is the contraction defined by
\begin{equation}\label{Ydef8}
Y = P_{_{\mathfrak{H}(\widetilde{V})}}\vert \mathfrak{H}(\widetilde{W})
: \mathfrak{H}(\widetilde{W}) \rightarrow \mathfrak{H}(\widetilde{V}).
\end{equation}
Recall that Proposition \ref{prop-R} shows that the dimensions of $\kr(Y)$ and
$ \kr (T_R) $ are equal.
Therefore we are interested in calculating the dimension of $\kr(Y)$ in terms
of the realizations of $V$ and $W$.
Since $H_W = \Gamma_{w} \Upsilon_w$, and the range of the  Hankel operator
$H_W^* $ is given by $\im(H_W^*) = \mathfrak{H}(\widetilde{W})$, we have
\[
\mathfrak{H}(\widetilde{W} ) =  \im ( \Upsilon_w^*)
\quad \mbox{where}\quad
 \Upsilon_w^* =
   B_w^* e^{A_w^* t} : \mathcal{X}_w \rightarrow L_+^2(\mathcal{E}).
\]
Furthermore, $\Upsilon_w $ is a co-isometry and $ A_w^* $
 is a stable dissipative operator on $\mathcal{X}_w$
satisfying $A_w^* +A_w + B_w B_w^* =0$, where $B_w$ is an operator mapping
$\mathcal{E}$ into $\mathcal{X}_w$.
Using the state space realization  $ V(s) = D_v +C_v(sI-A_v)^{-1}B_v$, we have
\begin{align*}
\left(T_{\widetilde{V}}^* \Upsilon_w^*\right)(t) &= \left(D_v \Upsilon_w^*\right)(t)
 + \int_t^\infty C_v e^{A_v (\tau-t)}B_v B_w^* e^{A_w^* \tau} d \tau\\
 &=\left(D_v \Upsilon_w^*\right)(t)
 +\left( \int_0^\infty C_v e^{A_v  \tau }B_v B_w^* e^{A_w^* \tau} d \tau \right)e^{A_w^* t}\\
 &=\left(D_v B_w^*
 +  C_v \int_0^\infty e^{A_v  \tau }B_v B_w^* e^{A_w^* \tau} d \tau\right) e^{A_w^* t}.
\end{align*}
In other words,
\begin{equation}\label{c000}
\left(T_{\widetilde{V}}^* \Upsilon_w^*\right)(t) =\left(D_v B_w^*
 +  C_v \int_0^\infty e^{A_v  \tau }B_v B_w^* e^{A_w^* \tau} d \tau\right) e^{A_w^* t}.
\end{equation}
Let $C_\circ$ be the operator mapping $\mathcal{X}_w$ into $\mathcal{E}$ defined by
\begin{equation}\label{defCo}
C_\circ
= D_v B_w^* + C_v\int_0^\infty e^{A_v  \tau }B_v B_w^* e^{A_w^* \tau} d \tau.
\end{equation}
Because the operators $A_v$ and $A_w$ are both stable,
the operator $C_\circ$ is well defined.
In fact, $C_\circ$ can be computed by first solving the Lyapunov
equation
\begin{equation}\label{lyapom}
 A_v \Omega + \Omega A_w^*  + B_v B_w^* =0.
\end{equation}
Because $A_v$ and $A_w$ are stable, the solution $\Omega$ to this
Lyapunov equation is unique and given by
\begin{equation} \label{defOmega}
\Omega = \int_0^\infty e^{A_v  \tau }B_v B_w^* e^{A_w^* \tau} d \tau.
\end{equation}
Therefore
\begin{equation}\label{Ccirc}
C_\circ = D_v B_w^* + C_v \Omega.
\end{equation}

Let $\Gamma_\circ$ be the operator
from $\mathcal{X}_w$ into $L_+^2(\mathcal{E})$
defined by  $ \Gamma_\circ = T_{\widetilde{V}}^* \Upsilon_w^* $.
By consulting \eqref{c000} with the definition of $C_\circ$ in \eqref{defCo},
we see that
\begin{equation}\label{Gammacirc}
 \left(\Gamma_\circ x\right)(t) =    \left(T_{\widetilde{V}}^* \Upsilon_w^* x\right)(t)  =
  C_\circ e^{A_w^* t}x \qquad (x\in \mathcal{X}_w).
\end{equation}
Notice  that $\Gamma_\circ$ is the observability operator
determined by the pair $\{C_\circ,A_w^*\}$.

It is emphasized that if $V = \varphi I$ where $\varphi$
is an inner function in $H^\infty$,
then $C_\circ = B_w^* \varphi(-A_w^*)$. In other words, in the scalar case
\begin{equation}\label{qqq}
\Gamma_\circ =
T_{\widetilde{\varphi}I}^* \Upsilon_w^* = \Upsilon_w^* \varphi(-A_w^*)
\end{equation}
and the operator $C_\circ$ plays the role of   $B_w^*\varphi(-A_w^*)$
in the MIMO case. Note that in the general case,
$C_\circ$ is the left point evaluation
of $-A_w^* $ with respect to $  B_w^* $.
For a further discussion on MIMO
function evaluation with applications to $H^\infty$ interpolation theory,
see Section 1.2 page 15 of \cite{FFGK}.

The next lemma  with $\mathfrak{n}(T_R) = \mathfrak{n}(Y)$
(see Part 2 of Proposition \ref{prop-R}) provides the proof of Part 1 and the first half of Part 2
of Theorem \ref{mainthm00}.

\begin{lemma}
Let $V$ and $W$ be  bi-inner rational  functions  in
$H^\infty(\mathcal{E},\mathcal{E})$.
Consider the contraction $Y$  mapping   $\mathfrak{H}(\widetilde{W})$
into $\mathfrak{H}(\widetilde{V})$  determined  by
\begin{equation}\label{Y43}
 Y = P_{_{\mathfrak{H}(\widetilde{V})}}\vert\mathfrak{H}(\widetilde{W}):
 \mathfrak{H}(\widetilde{W}) \rightarrow \mathfrak{H}(\widetilde{V}).
\end{equation}
Let $ C_\circ$ be the operator from $\mathcal{X}_w$ into $\mathcal{E}$
defined by \eqref{defCo} or \eqref{lyapom} and \eqref{Ccirc}.
Let $ Q$ be the unique solution to the Lyapunov equation
\begin{equation}\label{defQ}
A_w Q + Q A_w^* + C_\circ^* C_\circ =0.
\end{equation}
Then $Q$ is a positive contraction. Moreover, $Y$ and $I-Q$ have the same rank.
In particular,
\begin{equation}\label{eig1}
\kr(Y) = \Upsilon_w^*\kr(I - Q)
\quad \mbox{and}\quad \dim(\kr(Y))= \dim(\kr( I- Q)).
\end{equation}
Furthermore, we have
\begin{equation}\label{eig1dual}
\dim(\kr(Y^*))   = \dim(\mathfrak{H}(\widetilde{V})) - {\rm rank} (I -  Q).
\end{equation}
 If   $V = \varphi I$ where $\varphi$ is a Blaschke product,
 then $Q = \varphi(-A_w^*)^*\varphi(-A_w^*)$.
\end{lemma}

\bpr
Notice  that $ Q =  \Gamma_\circ^*  \Gamma_\circ $ is
the unique solution to the Lyapunov equation
\[
A_w Q +  Q A_w^* +  C_\circ^* C_\circ =0.
\]
Recall that the orthogonal projection
  $ P_{_{\mathfrak{H}(\widetilde{V})}}
= I - T_{\widetilde{V}} T_{\widetilde{V}}^* $.
Moreover,  $\mathfrak{H}(\widetilde{W})$ equals the
range of the isometry $\Upsilon_w^*$.
So $h$ is in $\mathfrak{H}(\widetilde{W})$
if and only if $h =  \Upsilon_w^* x$ for some
$x$ in  $\mathcal{X}_w$. In fact, $x = \Upsilon_w h$. Now observe that
\begin{align*}
  \|P_{_{\mathfrak{H}(\widetilde{V})}} \Upsilon_w^* x \|^2 &=
  \|\big(I - T_{\widetilde{V}} T_{\widetilde{V}}^* \big) \Upsilon_w^* x\|^2
  =
  \| \Upsilon_w^* x\|^2 -
  \| T_{\widetilde{V}} T_{\widetilde{V}}^* \Upsilon_w^* x \|^2\\
&=
\| \Upsilon_w^* x\|^2 - \|  T_{\widetilde{V}}^*\Upsilon_w^* x\|^2
 = \|x\|^2 - \|  \Gamma_\circ  x\|^2 = \|x\|^2 - \langle Q x,x\rangle.
\end{align*}
(The fourth equality follows from  \eqref{Gammacirc}.)
Hence $Q$ is a positive contraction, and thus,
\[
 \|P_{_{\mathfrak{H}(\widetilde{V})}} \Upsilon_w^* x \|^2
  = \|(I -  Q)^\frac{1}{2}x\|^2 \qquad (x\in \mathcal{X}).
\]
So there exists a unitary operator $\Psi$ mapping
the range of $P_{\mathfrak{H}(\widetilde{V} )} \Upsilon_w^* $ onto
the range $(I -  Q)^\frac{1}{2}$ such that
\begin{equation}\label{defect6}
\Psi Y  \Upsilon_w^* x =
\Psi P_{\mathfrak{H}(\widetilde{V} )}  \Upsilon_w^* x
= \big(I -  Q\big)^\frac{1}{2} x
\qquad (\mbox{for } x\in  \mathcal{X}_w ).
\end{equation}
In particular, $Y$ and $I-Q$ have the same rank.
Recall that $\mathfrak{H}( \widetilde{W} )$ equals the range of
$\Upsilon_w^*$
and that $Y = P_{_{\mathfrak{H}(\widetilde{V} )}} \vert  \mathfrak{H}(\widetilde{W}) $.
Therefore, by \eqref{defect6}
\[
\kr(Y) = \Upsilon_w^*\kr(I - Q)
\quad \mbox{and}\quad
\dim(\kr(Y))= \dim (\kr( I- Q)) .
\]
In other words, $\dim (\kr(Y) )$ equals the number of eigenvalues of $ Q$ equal to $1$ counting multiplicities.

By applying the previous result to $Y^*$, we obtain
\[
\dim (\kr(Y^*)) = \dim (\mathfrak{H}( \widetilde{V} ))  - \rank(I -  Q).
\]

Finally, it is noted that if $V = \varphi I$ where
$\varphi$ is an inner function in $H^\infty$, then
\[
\Gamma_\circ = \Upsilon_w^*\varphi(-A_w^*);
\]
see \eqref{qqq}. Because $\Upsilon_w $ is a co-isometry, we have
\[
Q = \Gamma_\circ^*\Gamma_\circ =
\varphi(-A_w^*)^* \Upsilon_w \Upsilon_w^* \varphi(-A_w^*)
= \varphi(-A_w^*)^* \varphi(-A_w^*).
\]
\epr

The following lemma implies Part 2 of Theorem \ref{mainthm00}.

\begin{lemma}
Let $V$ and $W$ be two rational bi-inner functions in
$H^\infty(\mathcal{E},\mathcal{E})$, and  $\psi$ be an
inner function in $H^\infty$.
Let $Y_\psi$ be the
contraction from $\mathfrak{H}(\widetilde{W})$ into
$\mathfrak{H}(\psi \widetilde{V})$ defined by
$Y_\psi = P_{\mathfrak{H}( \psi  \widetilde{V} )}
\vert {\mathfrak{H}( \widetilde{W} ) }$.
Then $\dim (\kr(Y_\psi )) $ is equal to the multiplicity of $1$ as an eigenvalue of
$\psi(-A_w) Q \psi(-A_w)^* $.
Also
\[
\dim (\kr(Y_\psi^*)) = \deg(\psi) \dim \mathcal{E} + \dim (\mathfrak{H}( \widetilde{V} )) -
\rank (I - \psi(-A_w) Q \psi(-A_w)^* ).
\]
In particular, if $\zeta(s) = \frac{1-s}{1+s}$, then
 $\dim (\kr(Y_{\zeta^k})) $ is equal to the multiplicity of $1$ as an eigenvalue of
$ \zeta(-A_w)^k Q \zeta(-A_w)^{*k} $.
Also
\[
\dim (\kr(Y_{\zeta^k}^*))  = k\dim (\mathcal{E} )+ \dim (\mathfrak{H}( \widetilde{V} )) -
\rank (I -  \zeta(-A_w)^k Q  \zeta(-A_w)^{*k} ).
\]
\end{lemma}

\bpr
Set $\Gamma = \Upsilon_w^*$.
Then for $x$ in $\mathcal{X}_w$ we have
\begin{equation*}
\Vert P_{\mathfrak{H}( \psi \widetilde{V}) }  \Gamma x \|^2 =
\| ( I - T_{  \psi\widetilde{ V}} T_{ \psi \widetilde{  V}}^* ) \Gamma x  \|^2 =
\| \Gamma x \|^2 - \| T_{\psi \widetilde{V}} T_{\psi \widetilde{V}}^* \Gamma x\|^2 .
\end{equation*}
Since  $ T_{ \psi\widetilde{ V}} $ is an isometry, and using \eqref{Gammacirc}, we have
\begin{align*}
\Vert P_{\mathfrak{H}(  \psi\widetilde{ V}) }  \Gamma x \|^2 &=
\| \Gamma x \|^2 - \| T_{ \psi \widetilde{  V}}^* \Gamma x\|^2
 = \| \Gamma x \|^2 - \| T_{  \widetilde{ V}}^* T_{\psi }^*\Gamma x\|^2\\
&= \| \Gamma x \|^2 - \| T_{  \widetilde{ V}}^*\Gamma \psi(-A_w)^*x\|^2\\
&= \| \Gamma x \|^2 - \| \Gamma_\circ \psi(-A_w)^*x\|^2.
\end{align*}
Recall that $ \Gamma_\circ^*  \Gamma_\circ  =  Q $.
Therefore
\begin{align*}
\Vert P_{\mathfrak{H}( \psi \widetilde{V}) }  \Gamma x \|^2 &=
\|x\|^2 - \langle \psi(-A_w) Q \psi(-A_w)^* x,x\rangle\\
&=
\|(I - \psi(-A_w) Q \psi(-A_w)^*)^\frac{1}{2}x\|^2.
\end{align*}
So there exists a unitary operator $\Psi$ mapping
the range of $P_{\mathfrak{H}( \psi \widetilde{V}) }  \Gamma $ onto
the range of
$(I - \psi(-A_w) Q \psi(-A_w)^*)^\frac{1}{2} $ such that
\[
\Psi P_{\mathfrak{H}( \psi \widetilde{V}) }  \Gamma x =
\big(I - \psi(-A_w) Q \psi(-A_w)^*\big)^\frac{1}{2} x
\qquad (\mbox{for } x \in  \mathcal{X}_w ).
\]
Recall that $\mathfrak{H}( \widetilde{W} ) = \Gamma \mathcal{X}_w$.
Therefore, when
$Y_\psi = P_{{\mathfrak{H}(\psi \widetilde{V} )}} \vert {\mathfrak{H}( \widetilde{W} )}$, we have
\[
\dim \left( \kr \left( P_{\mathfrak{H}(\psi \widetilde{V} )}
\vert {\mathfrak{H}( \widetilde{W})} \right)\right)
= \dim\left( \kr \left( I- \psi(-A_w) Q \psi(-A_w)^* \right) \right).
\]
In other words, the dimension of $ \kr(Y_\psi )$ is equal to the multiplicity of $1$ as an eigenvalue
of $ \psi(-A_w) Q \psi(-A_w)^* $.

Using $\mathfrak{H}(\psi \widetilde{V} )  = \mathfrak{H}(\psi I_\mathcal{E})\oplus
T_\psi \mathfrak{H}( \widetilde{V})$,
we also have
\[
\dim \left(\kr(Y_\psi^* ) \right) = \deg(\psi) \dim( \mathcal{E}) +
\dim (\mathfrak{H}( \widetilde{V} )) -
\mbox{rank}(I - \psi(-A_w) Q \psi(-A_w)^* ).
\]
\epr

 The previous lemma with  $\mathfrak{n}(T_{\zeta^kR}) = \mathfrak{n}(Y_{\zeta^k})$, yields
\[
\mathfrak{n}(T_{\zeta^k R})  = \dim \left( \kr (I - \zeta(-A_w)^{k} Q \zeta(-A_w)^{*k} )\right).
\]
This completes the proof of Part 2 of Theorem \ref{mainthm00}.

Since Part 3 of Theorem \ref{mainthm00} follows from the Parts 1 and 2 the proof of Theorem \ref{mainthm00} is complete now.

The next Corollary provides similar formulas for the positive Wiener-Hopf indices.
First we define $C_{\circ\ast}$ as follows.
Solve the Lyapunov
equation
\begin{equation}\label{lyapom_ast}
 A_w \Omega_\ast + \Omega_\ast A_v^* + B_w B_v^* =0.
\end{equation}
(Note that $\Omega_*=\Omega^*$.)
Because $A_v$ and $A_w$ are stable the solution $\Omega_*$ to this
Lyapunov equation is unique and given by
\begin{equation}
\Omega_\ast = \int_0^\infty   e^{A_w t} B_w B_v^* e^{A_vt} dt.
\end{equation}
Put
\begin{equation}\label{Ccircd}
C_{\circ\ast} = D_w B_v^* + C_w \Omega_\ast.
\end{equation}
Now let $ Q_\ast $ be the unique solution of the Lyapunov equation
\begin{equation}\label{defQast}
A_v Q_\ast + Q_\ast A_v^\ast + C_{\circ\ast}^* C_{\circ\ast} =0.
\end{equation}

\begin{corollary}\label{maincol}
Let $V$ and $ W $ be given by \eqref{defV00} and \eqref{defW00} and $ R = V W^* $.
Furthermore let $ \zeta(s) = \frac{1-s}{1+s}$ and $ Q_\ast $ be defined by \eqref{defQast}.
Put
\begin{align}\label{defnuk}
\nu_k & = \dim (\kr(I - \zeta(-A_v)^{k-1} Q_\ast \zeta(-A_v)^{*(k-1)} )) + \nonumber \\
& \quad - \dim ( \kr (I - \zeta(-A_v)^{k} Q_\ast \zeta(-A_v)^{*k} ) ).
\end{align}
Then the positive Wiener-Hopf indices $ \omega_1 , \ldots , \omega_q $
of $ T_R $ are given by
\[
\omega_j = \# \{ k \mid \nu_k \geq j \}, \quad (j = 1 , \ldots, q = \nu_1 ).
\]
\end{corollary}

\bpr
Notice $ R^* = W V^* $.
Recall $ R^*(s) $ is defined by $R^*(s) = \bigl(R(-\bar{s})\bigr)^* $.
If
\[
R(s) = W_-(s) \ \diag( \zeta(s)^{\kappa_j} )_{j=1}^m \ W_+(s) ,
\]
with  $ W_+ $ and its inverse are analytic on the right half plane $\BC_+$ and
$W_- $ and its inverse are analytic on the left half plane, $\BC_-$,
 then
\[
 R^*(s) = W_+^*(s) \diag( \zeta(s)^{-\kappa_j} )_{j=1}^m W_-^*(s) .
\]
Moreover $ W_-^* $ and its inverse are analytic on $ \BC_+ $ and
$W_+^* $ and its inverse are analytic on $\BC_- $.
This shows that $  - \kappa_1, \ldots -\kappa_m $
are the Wiener-Hopf indices of $ R^* $.
The positive Wiener-Hopf indices of $ R $ are the opposite to the
negative Wiener-Hopf indices of $ R^* $.
So the Corollary is immediate from applying Theorem \ref{mainthm00} to $ R^* $.
\epr

\bigskip\noindent

\setcounter{equation}{0}

\section{Rational matrix functions and the Cayley transform}
We denote the open unit disk in the complex plane by $\BD$,
the unit circle by $ \BT $, the closed unit disk by $ \overline{\BD} $,
the open right hand half plane by $ \BC_+ $, the imaginary axis by $ i \BR $
and the open left hand half plane by $ \BC_- $.

Let the conformal transformation $ \zeta : \BC \to \BC $ be given by
\begin{equation*}
z = \zeta( s ) = \textstyle{ \frac{1-s }{s+1 } }  \quad ( s\in \BC \setminus \{-1\}).
\end{equation*}
Then also
\[
s = \zeta^{-1} (z) = \textstyle{ \frac{1-z }{z+1 } }\quad ( z\in \BC \setminus \{-1\}).
\]
If $ s \in \BC_+ $ then the distance of $ s $ to $ -1$
is larger than the distance to $1$,
or in formula $ | s+1 | > | 1-s| $ or $ |\zeta(s)| <1 $.
So $ \zeta (\BC_+ ) = \BD $.
For $ s \in i\BR $ the distance of $ s $ to $ -1$
is equal to the distance to $1$,
or in formula $ | s+1 | = | 1-s| $ or $ |\zeta(s)| = 1 $.
So $ \zeta( i\BR ) = \BT $.
Similarly one sees $ \zeta( \BC_- ) = \BC \setminus \overline{\BD} $.
Using the similarity of the formula of $ \zeta^{-1} $
with the formula for $\zeta$ we see that $ \zeta (\BD ) = \BC_+ $,
$ \zeta( \BT ) = i\BR $ and $ \zeta( \BC \setminus \overline{\BD} ) = \BC_- $.

The following  result  is well-known, compare Section 3.6 in \cite{BGKROT178}. The result
is presented here for sake of completeness.

\begin{lemma}\label{DtoC}
Let $\{A_d \mbox{ on } \mathcal{X},B_d,C_d,D_d \}$
be a discrete time realization of the  function $ \Theta $, that is,
 \[
 \Theta(z) = D_d + z C_d (I - z A_d)^{-1} B_d.
 \]
Then
\begin{equation}\label{thetasigma}
\Theta(\zeta(s)) = D_c + C_c ( s I - A_c)^{-1} B_c.
\end{equation}
Here
\begin{align}
 A_c &= (A_d -I)(I + A_d )^{-1}, \qquad
 D_c = D_d - C_d ( I + A_d )^{-1}B_d, \label{AD}\\
 B_c &= \sqrt{2}(I+A_d)^{-1}B_d,\qquad \quad \
 C_c = \sqrt{2} C_d (I+A_d )^{-1}.\label{BC}
\end{align}
\end{lemma}
\bpr
First note that
\begin{align*}
\Theta(\zeta(s)) &=
\textstyle{ D_d + \frac{1-s}{ s +1} C_d \left( I - \frac{1-s}{ s +1} A_d \right)^{-1}  B_d   }\\
& = D_d + (1-s) C_d \big(  (s+1 ) I - ( 1-s ) A_d \big)^{-1} B_d  \\
& = D_d + ( 1-s) C_d \big( (I - A_d ) + s (I + A_d ) \big)^{-1} B_d.
\end{align*}
By employing $ A_c = (A_d - I )(I + A_d)^{-1} $, we have
\begin{equation*}
\Theta(\zeta(s)) =
D_d + ( 1-s) C_d \left( sI - A_c \right)^{-1} (I+A_d)^{-1} B_d.
\end{equation*}
Since $ (1-s) \left( sI - A_c \right)^{-1} = (I - A_c )( sI - A_c)^{-1} -I$, we obtain
\begin{equation*}
\Theta(\zeta(s))= D_d - C_d (I+ A_d )^{-1} B_d + C_d (I - A_c )\left( sI -  A_c \right)^{-1} (I+A_d)^{-1} B_d.
\end{equation*}
Finally, by using $ I - A_c = 2 (I+A_d)^{-1} $, we see that
\begin{equation*}
\Theta(\zeta(s)) =
D_d -C_d  (I+A_d)^{-1} B_d +
 C_d \sqrt{2} (I + A_d )^{-1} \left( sI -A_c \right)^{-1} \sqrt{2}(I+A_d)^{-1} B_d.
\end{equation*}
We conclude that \eqref{thetasigma} holds.
\epr

\medskip

\medskip
We can also express the operators $A_d$, $B_d$, $C_d$ and $D_d$ in
$A_c$, $ B_c $, $ C_c $, $D_c$. Indeed we have
the following result.
\begin{corollary}\label{CtoD}
With $A_d$, $B_d$, $C_d$ $D_d$ and $A_c$, $ B_c $, $ C_c $, $D_c$ as in the Lemma~\ref{DtoC} we have
\begin{align*}
& A_d = (I - A_c )^{-1} (A_c + I ), \quad
B_d = \sqrt{2} \, (I - A_c)^{-1} B_c,  \\
& C_d = \sqrt{2} \, C_c ( I - A_c )^{-1} , \quad
D_d =D_c + C_c (I - A_c )^{-1} B_c.
\end{align*}
In particular, if $G(s)$ admits a continuous time realization of the form
\[
G(s) = D_c + C_c(sI-A_c)^{-1}B_c,
\]
then the corresponding discrete time realization is given by
\[
G(\zeta^{-1}(z)) = D_d +z C_d(I-z A_d)^{-1}B_d,
\]
and vice-versa.
\end{corollary}

\bpr
This follows from a simple rewriting of \eqref{AD} and \eqref{BC}.
\epr

\medskip

From \cite{sznf, sznfbk} it is known that inner functions on the unit circle
have a stable unitary realization, and conversely.
The Cayley transform converts a stable unitary realization into
a stable dissipative realization as is shown in the following lemma.

Recall that a discrete time transfer function $\Theta(z)$ is bi-inner if and only if
$\Theta(z)$ admits a stable discrete time realization
\begin{align}\label{unitary}
&\Theta(z) = D_d +z C_d(I-z A_d)^{-1}B_d \qquad \mbox{where}\quad \nonumber\\
&\begin{bmatrix} A_d & B_d \\ C_d & D_d  \end{bmatrix}:\begin{bmatrix}\mathcal{X}_d \\ \mathcal{U}  \end{bmatrix}
\rightarrow\begin{bmatrix}\mathcal{X}_d \\ \mathcal{U}  \end{bmatrix} \qquad \qquad\,\,\mbox{is unitary}.
\end{align}
Finally, a discrete time system $\{A_d,B_d,C_d,D_d\}$ is \emph{stable} if all the eigenvalues
for $A_d$ are contained in $\mathbb{D}$.

\begin{lemma}\label{diss-sta}
Let $\{A_d \mbox{ on } \mathcal{X},B_d,C_d,D_d \}$
be a discrete time  realization of the  function $ \Theta(z) $, that is,
 \[
 \Theta(z) = D_d + z C_d (I - z A_d)^{-1} B_d.
 \]
Let  $\{A_c \mbox{ on } \mathcal{X},B_c,C_c,D_c \}$ be the corresponding
continuous time realization of $ \Theta( \zeta(s)) $,  where
\begin{align*}
 A_c &= (A_d -I)(I + A_d )^{-1}, \quad \quad
 D_c = D_d - C_d ( I + A_d )^{-1}B_d,\\
 B_c &= \sqrt{2}(I+A_d)^{-1}B_d, \qquad \,
\quad  C_c = \sqrt{2} C_d (I+A_d )^{-1}.
\end{align*}
Then $\{A_d \mbox{ on } \mathcal{X},B_d,C_d,D_d \}$ is
a stable unitary realization of $ \Theta(z) $ if and only if
$\{A_c \mbox{ on } \mathcal{X},B_c,C_c,D_c \}$
is a stable dissipative realization of $ \Theta(\zeta(s)) $.
\end{lemma}

\bpr
The continuous time realization $\{A_c \mbox{ on } \mathcal{X},B_c,C_c,D_c \}$ being
stable and dissipative means that $ A_c $ has all its eigenvalues
in the left hand half plane $ \BC_- $ and $ A_c + A_c^\ast + C_c^* C_c =0 $,
$D_c$ is unitary, and $ B_c = - C_c^\ast D_c$.

To start with the stability  we suppose that $ A_c x = \lambda x $.
Then
\[
A_d x = (A_c -I)^{-1} ( A_c +I)x = \frac{ 1+ \lambda }{\lambda -1} x.
\]
Note that $ \frac{ 1+ \lambda }{\lambda -1}$ is an eigenvalue of $A_d $
and this means $\frac{ 1+ \lambda }{\lambda -1} \in \BD $.
In other words $ | \lambda +1 | < |\lambda -1| $ and this is equivalent to
 $ \lambda \in \BC_- $.
We conclude that $ A_d $ is (discrete time) stable if and only if $ A_c $ is (continuous time) stable.

Recall that \eqref{unitary} being unitary means that
\begin{equation}\label{d-unitary}
\begin{bmatrix}   I & 0 \\  0 & I \\\end{bmatrix}
= \begin{bmatrix}  A_d^* & C_d^* \\   B_d^* & D_d^* \\ \end{bmatrix}
\begin{bmatrix}   A_d & B_d \\  C_d & D_d \\\end{bmatrix}
= \begin{bmatrix} A_d^* A_d + C_d^*C_d & A_d^*B_d +C_d^*D_d \\
       B_d^*A_d + D_d^* C_d & B_d^*B_d + D_d^*D_d \end{bmatrix}.
\end{equation}

Notice that
\begin{equation}\label{CcCc}
C_c^* C_c = 2 (I+A_d^*)^{-1} C_d^* C_d (I + A_d )^{-1}
\end{equation}
and
\begin{align}
A_c + A_c^* &= (A_d - I) ( I+A_d)^{-1} + ( I+A_d^* )^{-1} ( A_d^* -I ) \nonumber\\
&= ( I+A_d^* )^{-1}
\bigl( (A_d^* - I)  ( I+A_d) +( I+A_d^* )(A_d - I) \bigr)
 ( I+A_d)^{-1} \nonumber\\
 &= (I+A_d^*)^{-1} ( 2 A_d^* A_d -2I )(I+ A_d )^{-1} . \label{Ac+Ac}
\end{align}
Adding the two equalities \eqref{CcCc} and \eqref{Ac+Ac} we see that
\begin{equation*}
 A_c + A_c^\ast + C_c^* C_c =
2(I+A_d^*)^{-1} \bigl(  A_d^* A_d + C_d^* C_d -I \bigr)(I+ A_d )^{-1} .
\end{equation*}
We conclude that
$ A_c + A_c^\ast + C_c^* C_c =0 $ if and only if $ A_d^* A_d + C_d^* C_d -I =0 $.

The following calculation shows that $ B_c = - C_c^\ast D_c$ if and only if\\
$ A_d^* B_d + C_d^* D_d = 0 $:
\begin{align*}
& A_d^* B_d + C_d^* D_d =
\sqrt{2} \Big[ (I-A_c^*)^{-1} (I+ A_c^*) (I - A_c)^{-1} B_c + (I-A_c^*)^{-1} C_c^* D_c +\\
&\qquad + (I-A_c^*)^{-1}C_c^*C_c(I-A_c)^{-1} B_c   \Big] \\
&= \sqrt{2} (I-A_c^*)^{-1}  \big[ (I+ A_c^*) (I - A_c)^{-1} B_c + C_c^*C_c(I-A_c)^{-1} B_c +  C_c^* D_c  \big]  \\
&= \sqrt{2} (I-A_c^*)^{-1}
\Big[ (I+ A_c^*) (I - A_c)^{-1} B_c -(A_c + A_c^*) (I-A_c)^{-1} B_c +  C_c^* D_c  \big] \\
&=  \sqrt{2} (I-A_c^*)^{-1}  ( B_c + C_c^* D_c ).
\end{align*}

Finally, we check that, given $ A_c + A_c^\ast + C_c^* C_c =0 $ and
$ B_c = - C_c^\ast D_c$, we have that $ D_d^* D_d + B_d^* B_d = I $ if and only if
$ D_c^* D_c = I  $, i.e., $ D_c $ is unitary.
To this end, observe that
\begin{align*}
& D_d^* D_d + B_d^* B_d = \bigl[ D_c^* + B^* (I-A_c^*)^{-1}C_c^*  \bigr]
\bigl[ D_c + C_c (I-A_c^*)^{-1} B_c  \bigr] + \\
&\qquad+ 2 B_c^* (I-A_c^*)^{-1} (I-A_c)^{-1} B_c \\
&= D_c^* D_c + D_c^* C_c (I-A_c)^{-1} B_c + B_c^* (I-A_c^*)^{-1}C_c^* D_c \\
&\qquad + B_c^* (I-A_c^*)^{-1}  C_c^* C_c  (I-A_c)^{-1} B_c +
2 B_c^* (I-A_c^*)^{-1} (I-A_c)^{-1} B_c \\
&=  D_c^* D_c  + B_c^* (I-A_c^*)^{-1}
\Bigl[ -(I-A_c^*) - (I-A_c )- A_c - A_c^* + 2I  \Bigr] (I-A_c)^{-1} B_c \\
&=  D_c^* D_c ,
\end{align*}
In the third  equality we used that $ B_c = - C_c^* D_c $ and
$ A_c + A_c^* + C_c^* C_c =0 $.
Since $ D_d^* D_d + B_d^* B_d  = D_c^* D_c  $ we conclude that
$ D_d^* D_d + B_d^* B_d  =I  $ if and only if $  D_c^* D_c  = I  $.

Therefore $\{A_d \mbox{ on } \mathcal{X},B_d,C_d,D_d \}$ is
a (discrete time) stable unitary realization of $ \Theta(z) $ if and only if
$\{A_c \mbox{ on } \mathcal{X},B_c,C_c,D_c \}$
is a stable dissipative (continuous time) realization of $ \Theta(\zeta(s)) $.
\epr

\bigskip\bigskip

\section{Equivalence of discrete time and continuous time theorems}
Consider   the functions   $ V(z) $ and  $W(z)$ in $H_{\mathbb{D}}^\infty(\mathcal{E},\mathcal{E})$.
Let  $V(\zeta(s)) $ and $W(\zeta(s) )$ in $H^\infty(\mathcal{E},\mathcal{E})$ be their
corresponding Cayley transforms.   Moreover, consider their corresponding discrete and
continuous time realizations given by
\begin{align}
 V(z) &= D_{dv} + z C_{dv} (I - z A_{dv})^{-1} B_{dv}, \label{reaVd} \\
  W(z) &= D_{dw} + z C_{dw} (I - z A_{dw})^{-1} B_{dw}, \label{reaWd} \\
 V(\zeta(s)) &= D_{cv} + C_{cv}  ( s I - A_{cv} )^{-1} B_{cv}, \label{reaVc} \\
 W(\zeta(s)) &= D_{cw} + C_{cw} ( s I - A_{cw})^{-1} B_{cw}.  \label{reaWc}
\end{align}
The discrete realization  in \eqref{reaVd}  and  \eqref{reaWd},  and their
corresponding continuous time counter parts  \eqref{reaVc} and \eqref{reaWc}
are related  by the transformations presented in  Lemma \ref{DtoC}.
Recall from Corollary \ref{CtoD}
\begin{align}\label{need1}
& A_{dv} = (I - A_{cv} )^{-1} (A_{cv} + I ), \quad
B_{dv} = \sqrt{2} \, (I - A_{cv})^{-1} B_{cv},  \nonumber \\
& C_{dv} = \sqrt{2} \, C_{cv} ( I - A_{cv} )^{-1} , \quad
D_{dv} =D_{cv} + C_{cv} (I - A_{cv} )^{-1} B_{cv},
\end{align}
and
\begin{align}\label{need2}
& A_{dw} = (I - A_{cw} )^{-1} (A_{cw} + I ), \quad
B_{dw} = \sqrt{2} \, (I - A_{cw})^{-1} B_{cw}, \nonumber \\
& C_{dw} = \sqrt{2} \, C_{cw} ( I - A_{cw} )^{-1} , \quad
D_{dw} =D_{cw} + C_{cw} (I - A_{cw} )^{-1} B_{cw}.
\end{align}

Next we would like to develop a relationship between
a special Stein equation in the discrete time and a
corresponding Lyapunov equation in continuous time. To this end,
consider the Stein equation:
\begin{equation}\label{OmegaD}
\Omega = A_{dv} \Omega A_{dw}^\ast + B_{dv} B_{dw}^\ast.
\end{equation}
By employing the corresponding transformations in
\eqref{need1} and \eqref{need2}, we obtain
\begin{align*}
\Omega & =
(I - A_{cv} )^{-1} (A_{cv} + I ) \Omega (A_{cw}^\ast + I )(I - A_{cw}^\ast )^{-1}
 + \\
& \quad + (I - A_{cv})^{-1} B_{cv} 2 B_{cw}^\ast  (I - A_{cw}^\ast)^{-1} .
\end{align*}
Multiplying   by $ I - A_{cv} $ on the left and  by $ I - A_{cw}^\ast $ on the right, we
arrive at
\[
(I - A_{cv}) \Omega (I - A_{cw}^\ast)  =
 (A_{cv} + I ) \Omega (A_{cw}^\ast + I ) + 2 B_{cv}  B_{cw}^\ast .
\]
This simplifies to
\begin{equation} \label{OmegaC}
A_{cv} \Omega + \Omega A_{cw}^\ast + B_{cv}  B_{cw}^\ast = 0 .
\end{equation}
Therefore  $ \Omega $ is a solution of the Stein equation  \eqref{OmegaD} if and only if
$ \Omega $  is a solution of the Lyapunov equation \eqref{OmegaC}.

Next we establish the relation between $ C_{c\circ} $ and $ C_{d\circ} $ where
\[
C_{d\circ} = D_{dv} B_{dw}^\ast + C_{dv} \Omega A_{dw}^\ast
\]
and
\[
C_{c\circ} = D_{cv} B_{cw}^\ast + C_{cv} \Omega.
\]
We have
\begin{align*}
C_{d\circ} & = D_{dv} B_{dw}^\ast + C_{dv} \Omega A_{dw}^\ast \\
&= \sqrt{2} \Bigl[  [ D_{cv} + C_{cv} (I - A_{cv} )^{-1} B_{cv} ] B_{cw}^\ast + \\
&\qquad + C_{cv} (I - A_{cv} )^{-1} \Omega (I + A_{cw}^\ast ) \Bigr] (I - A_{cw}^\ast )^{-1} .
\end{align*}
Now use \eqref{OmegaC} in
\begin{align*}
& C_{cv} (I - A_{cv} )^{-1}  \bigl[ B_{cv} B_{cw}^\ast + \Omega (I + A_{cw}^\ast )
\bigr] \\
&\quad = C_{cv} (I - A_{cv} )^{-1} \bigl[  (I - A_{cv} ) \Omega \bigr]\\
&\quad = C_{cv} \Omega .
\end{align*}
So we have
\begin{align*}
C_{d\circ} & = \sqrt{2} \Bigl[ D_{cv} B_{cw}^\ast (I - A_{cw}^\ast )^{-1}   + C_{cv} \Omega (I - A_{cw}^\ast )^{-1} \Bigr] .
\end{align*}
We conclude that
\begin{equation}\label{CcirctoDcirc}
C_{d\circ} =  \sqrt{2} C_{c\circ} (I - A_{cw}^\ast )^{-1} .
\end{equation}

Our aim is to develop a connection between \cite[Theorem 2.2]{FKRvS}
in the discrete time setting and our Theorem \ref{mainthm00}.
To this end, consider the equations
\begin{equation}\label{Eq-Qd}
Q_d = A_{dw} Q_d  A_{dw}^\ast + C_{d\circ}^\ast C_{d\circ}
\end{equation}
and
\begin{equation}\label{Eq-Qc}
A_{cw} Q_c + Q_c A_{cw}^\ast + C_{c\circ}^\ast C_{c\circ} =0.
\end{equation}
Because $A_{dw}$ is discrete time stable and $A_{cw}$ is
continuous time stable, the solution to both of these
equations is unique. We claim that $Q_d =  Q_c $.
Consider
\begin{align*}
Q & = A_{dw} Q  A_{dw}^\ast + C_{d\circ}^\ast C_{d\circ} \\
&=(I - A_{cw} )^{-1} (A_{cw} + I ) Q (A_{cw}^\ast + I ) (I - A_{cw}^\ast )^{-1} \\
&\qquad +2 (I - A_{cw} )^{-1}C_{c\circ}^\ast C_{c\circ} (I - A_{cw}^\ast )^{-1} .
\end{align*}
Multiply on the left with $ I - A_{cw} $ and on the right with $ I - A_{cw}^\ast $,
we have
\begin{equation*}
(I - A_{cw} ) Q (I - A_{cw}^\ast )  = (A_{cw} + I ) Q (A_{cw}^\ast + I )
 + 2  C_{c\circ}^\ast C_{c\circ} .
\end{equation*}
This simplifies to
\begin{equation*}
 A_{cw}  Q + Q A_{cw}^\ast  + C_{c\circ}^\ast C_{c\circ} =0 .
\end{equation*}
We conclude that the equations \eqref{Eq-Qc} and \eqref{Eq-Qd} have the same solutions.

\bigskip

We will show the equivalence of the next two Theorems.
The first is a rephrase of \cite[Theorem 2.2]{FKRvS} and the second rephrases Theorem \ref{mainthm00} above.
We denote with $H^\infty_{\mathbb{D}} (\mathcal{E},\mathcal{E}) $ the space
consisting of the set of all operator valued functions
$\Theta(s)$ on $\mathcal{E}$ that are analytic in
the open unit circle $\BD$ and such that
\[
\|\Theta\|_\infty = \sup \{ \|\Theta(z) \|: z \in \BD \} <\infty.
\]
In the remainder of the paper the boldface  $\mathbf{T}_R$ denotes
the Toeplitz operator on $\ell_+^2(\mathcal{E})$ with symbol
 $R$ in $L_\mathbb{T}^\infty(\mathcal{E},\mathcal{E})$.

\begin{theorem}\label{mainthm00d}  
Assume that $R(e^{i\omega}) = V(e^{i\omega})  W(e^{i\omega}) ^*$ where
 $V $ and $W $ are two bi-inner rational functions in $H_{\mathbb{D}}^\infty(\mathcal{E},\mathcal{E})$. Let $\mathbf{T}_R$ be the corresponding Toeplitz operator on $\ell_+^2(\mathcal{E})$.
Let $\{A_{dv}, B_{dv},C_{dv},D_{dv} \}$ and
$\{A_{dw}, B_{dw},C_{dw},  D_{dw}\}$  be
stable unitary realizations of $V(z)$ and $W(z)$,  respectively.
Let $\Omega_d $ be the unique solution of the Stein equation
\begin{equation}\label{lyapom00d}
\Omega_d = A_{dv} \Omega_d A_{dw}^* + B_{dv} B_{dw}^*.
\end{equation}
Let $C_{d\circ}$ be the operator mapping $\mathcal{X}_w$ into $\mathcal{E}$ defined by
\begin{equation}
C_{d\circ} = D_{dv} B_{dw}^* + C_{dv} \Omega_d A_{dw}^*.
\end{equation}
Finally,  let $Q_d$ be the unique solution to the Stein equation
\begin{equation}\label{defQ00d}
 Q_d = A_{dw}  Q_d A_{dw}^* + C_{d\circ}^* C_{d\circ}.
\end{equation}
Then the following holds:
\begin{enumerate}
\item The operator $Q_d$ is a positive contraction.
\item The multiplicity of $1$ as an eigenvalue of $ Q_d $ equals $\mathfrak{n}(\mathbf{T}_R)$.
In other words, $\mathfrak{n}(T_R) = \mathfrak{n}(I-Q)$.
Moreover, for $k=0,1,2,\cdots, p$, we have
\begin{equation}\label{alphathm00d}
\mathfrak{n}(\mathbf{T}_{z^k R})  = \dim \left( \kr (I - A_{dw}^{k} Q A_{dw}^{*k} )\right).
\end{equation}
 \item For $k=1,2,\cdots$, consider the integers
\begin{equation}\label{defmuk00d}
\mu_k = \mathfrak{n}(I - A_{dw}^{k-1} Q (A_{dw}^*)^{k-1} )
- \mathfrak{n} (I - A_{dw}^{k} Q A_{dw}^{*k} ) .
\end{equation}
Then the negative Wiener-Hopf indices $ -\kappa_1 , \ldots , -\kappa_p $
of the Toeplitz operator $\mathbf{T}_R $ are given by
\begin{equation}\label{defmuk001d}
\kappa_j = \# \{ k : \mu_k \geq j \}, \quad (j = 1 , \ldots, p = \mu_1 ).
\end{equation}
\end{enumerate}
\end{theorem}

\bigskip

Recall that the  transformation $ \zeta : \BC \to \BC $ is given by
\begin{equation*}
\zeta(s) = \textstyle{ \frac{1-s }{s+1 } } \quad ( s\in \BC \setminus \{-1\}).
\end{equation*}

\begin{theorem}\label{mainthm00c} 
Assume that $R(\zeta(i\omega)) = V(\zeta(i\omega)) W(\zeta(i\omega))^*$
( $ \omega \in \BR $ ) where
 $V(\zeta(s))$ and $W(\zeta(s))$ are two bi-inner rational functions in $H^\infty(\mathcal{E},\mathcal{E})$.
Let
\[
\{A_{cv}, B_{cv},C_{cv},D_{cv} \}  \quad \mbox{and}\quad
\{A_{cw}, B_{cw},C_{cw},  D_{cw}\}
\]
 be two
stable dissipative  realizations of $V(\zeta(s))$ and $W(\zeta(s))$,  respectively.
Let $\Omega_c $ be the unique solution of the Lyapunov  equation
\begin{equation}\label{lyapom00c}
  A_{cv} \Omega_c +\Omega_c A_{cw}^* + B_{cv} B_{cw}^*=0.
\end{equation}
Let $C_{c\circ}$ be the operator mapping $\mathcal{X}_w$ into $\mathcal{E}$ defined by
\begin{equation}
C_{c\circ} = D_{cv} B_{cw}^* + C_{cv} \Omega_c.
\end{equation}
Finally,  let $Q_c $ be the unique solution to the Lyapunov  equation
\begin{equation}\label{defQ00c}
 A_{cw}  Q_c +Q_c  A_{cw}^* + C_{c\circ}^* C_{c\circ} =0.
\end{equation}
Then the following holds:
\begin{enumerate}
\item The operator $Q_c$ is a positive contraction.
\item The multiplicity of $1$ as an eigenvalue of $ Q_c $ equals
$\mathfrak{n}(T_{R\zeta})$.
In other words, $\mathfrak{n}(T_{R\zeta}) = \mathfrak{n}(I-Q_c)$.
Moreover, for $\zeta(s) = \frac{1-s}{1+s}$, we have
\begin{equation}\label{alphathm00c}
\mathfrak{n}(T_{\zeta^k (R\zeta)})  = \dim \left( \kr (I - \zeta(-A_{cw})^{k} Q_c \zeta(-A_{cw})^{*k} )\right).
\end{equation}
  \item For $k=1,2,\cdots$, consider the integers
\begin{align}\label{defmuk00c}
&\mu_k =
\mathfrak{n}(I - \zeta(-A_{cw})^{k-1} Q_c (\zeta(-A_{cw})^*)^{k-1} ) +  \nonumber\\
& \qquad \qquad - \mathfrak{n} (I - \zeta(-A_{cw})^{k} Q_c \zeta(-A_{cw})^{*k} ) .
\end{align}
Then the negative Wiener-Hopf indices $ -\kappa_1 , \ldots , -\kappa_p $
of $ T_{R\zeta} $ are given by
\begin{equation}\label{defmuk001c}
\kappa_j = \# \{ k : \mu_k \geq j \}, \quad (j = 1 , \ldots, p = \mu_1 ).
\end{equation}
\end{enumerate}
\end{theorem}

Assume we have Theoren \ref{mainthm00d}.
Since  $ Q_c = Q_d $, Part 1 of Theorem \ref{mainthm00c} is proven.
Notice that the Wiener-Hopf factorization of $ R$ with respect to $ \BT $
immediately  generates a Wiener-Hopf factorization of
$ R(\zeta(i\omega)) $ with respect to $ i\BR $ with the same Wiener-Hopf indices.
Indeed, if $ R = W_+ D W_- $ then,
after substitution of $z$ by $ \zeta(s) $, we have
$ R\zeta = (W_+\zeta) (D\zeta) (W_-\zeta)   $.
So the dimensions of the null spaces of $ \mathbf{T}_{z^k R} $
and of $ T_{\zeta^k R\zeta} $ coincide.
Together with $ A_{dw} = \zeta( -A_{cw}) $  the items 2 and 3
of Theorem \ref{mainthm00d} give the items 2 an 3 of Theorem \ref{mainthm00c}.

Obviously, also Theorem \ref{mainthm00c} implies Theorem \ref{mainthm00d}.

\subsection{An example}
As an illustration of Theorem \ref{mainthm00c} we present the following example
which is a continuous version of the discrete example on page 706 in \cite{GKRa}.
To this end, let
\[
R_c(s) = \begin{bmatrix}
\zeta(s)^{-4}&0&0&0&0\\  0& \zeta(s)^{-2}&0&0&0\\  0&0&1&0&0\\  0&0&0&\zeta(s)^3&0\\  0&0&0&0&\zeta(s)^5
\end{bmatrix}.
\]
Recall that $\zeta(s) = \frac{1-s}{1+s}$.
Then $ R_c(s) $ factors as $ R_c(i \omega) = V_c(i \omega) W_c(i\omega)^*$, where
\[
V_c(s) = \begin{bmatrix}
1&0&0&0&0\\  0&1&0&0&0\\  0&0&1&0&0\\  0&0&0&\zeta(s)^3&0\\  0&0&0&0&\zeta(s)^5
\end{bmatrix}
\qquad
W_c(s) = \begin{bmatrix}
\zeta(s)^{4}&0&0&0&0\\  0&\zeta(s)^{2}&0&0&0\\  0&0&1&0&0\\  0&0&0&1&0\\  0&0&0&0&1
\end{bmatrix} .
\]

Let us first construct a
stable dissipative realization $\{A_n \mbox{ on } \mathbb{C}^n,B_n,C_n,D_n\}$ for $\zeta(s)^n$.
Notice that $p(z) = z e_1^T \left(I -z J_n(0) \right)^{-1} e_n = z^n$,
where $e_1 $ and $e_n $ are the first and last unit vectors in $ \BC^n $
and $ J_n(0) $ is the upper triangular Jordan matrix  with eigenvalue 0,
or, in other words, the upward shift.
Then Lemma \ref{DtoC}  provides the realization of
$ \zeta(s)^n = p (\zeta(s))$.
To see this observe that
\[
 ( I_n+ J_n(0) )^{-1} = \sum_{j=0}^{n-1} ( - J_n(0) )^j.
 \]
 This with Lemma \ref{DtoC} readily implies that
\begin{align}\label{need}
A_n &= ( J_n(0) - I_n) ( I_n+ J_n(0) )^{-1}=
-I_n + 2 \sum_{j=1}^{n-1} (-1)^{j+1} J_n(0)^j , \nonumber\\
 B_n &= \sqrt{2} \,(I_n + J_n(0) )^{-1} e_n = \sqrt{2}
\left[ (-1)^{n-1} \ \ (-1)^{n-2} \ldots \ -1\ \ 1 \right]^T, \nonumber\\
 C_n &= \sqrt{2} \  e_1^T ( I_n+ J_n(0) )^{-1} =
\sqrt{2} \left[ 1 \ \ -1 \  \ldots \ (-1)^{n-1}  \right] , \nonumber\\
 D_n &= 0 - e_1^T (I_n + J_n(0) )^{-1} e_n = (-1)^n .
\end{align}
Notice that $ A_n $ is an upper triangular Toeplitz matrix.

According to Lemma \ref{diss-sta} and using that the realization
of $ p(z) $ is unitary, we have that
$\{A_n,B_n,C_n,D_n\}$ is indeed a
stable dissipative  realization of $\zeta(s)^n$.
Nevertheless let us verify that directly.
Recall that the realization $\{A_n \mbox{ on } \mathcal{X},B_n,C_n,D_n \}$ being
stable and dissipative means that $ A_n $ has all its eigenvalues
in the open left hand half plane  and $ A_n + A_n^\ast + C_n^* C_n =0 $,
$D_n$ is unitary, and $ B_n = - C_n^\ast D_n$.
First observe that $-1$ is the only eigenvalue of  $A_n $,
and thus, $A_n$ is continuous time stable.
Next check by direct calculation that $ A_n^* +A_n +C_n^*C_n = 0$ and
$A_n^* +A_n + B_n B_n^* = 0$. Since $D_n = (-1)^n$ is unitary and
$B_n = -C_n^* D_n$ we are done.

Motivated by the previous realization, the factors $ V_c(s) $ and $ W_c(s) $
can be given by the following stable unitary realizations:
\[
V_c(s) = D_v +   C_v ( s I  -   A_v )^{-1} B_v .
\]
where
\begin{align*}
 A_v &= \begin{bmatrix} A_3 & 0 \\ 0 & A_5 \end{bmatrix} :
\begin{bmatrix} \mathbb{C}^3   \\ \mathbb{C}^5 \end{bmatrix} \to
\begin{bmatrix} \mathbb{C}^3   \\ \mathbb{C}^5\end{bmatrix} , \\
B_v &= \begin{bmatrix} 0_{ 3 \times 3 } & B_3 & 0  \\
                       0_{ 5 \times 3 } & 0 &  B_5  \end{bmatrix}:
\begin{bmatrix}  \mathbb{C}^3 \\  \mathbb{C}  \\  \mathbb{C} \end{bmatrix}  \to \begin{bmatrix} \mathbb{C}^3   \\ \mathbb{C}^5 \end{bmatrix} , \\
C_v &=
\begin{bmatrix}0_{3 \times 3} &0_{3 \times 5} \\ C_3 & 0 \\ 0 & C_5 \end{bmatrix}:
 \begin{bmatrix}  \mathbb{C}^3 \\  \mathbb{C}^5 \\  \end{bmatrix}\rightarrow \begin{bmatrix}  \mathbb{C}^3 \\  \mathbb{C}  \\  \mathbb{C}  \end{bmatrix}
 = \BC^5 , \\
 D_v &=
\begin{bmatrix} I_3 & 0 & 0 \\ 0 & D_3 & 0 \\ 0 & 0 & D_5  \end{bmatrix}:
\begin{bmatrix} \mathbb{C}^3 \\ \mathbb{C}  \\ \mathbb{C} \end{bmatrix}\to   \begin{bmatrix} \mathbb{C}^3 \\ \mathbb{C}  \\ \mathbb{C}  \end{bmatrix}.
\end{align*}
In this case,
\[
W_c(s) = D_w +   C_w ( s I  -   A_w )^{-1} B_w ,
\]
where
\begin{align*}
 A_w &= \begin{bmatrix} A_4 & 0 \\ 0 & A_2 \end{bmatrix} :
\begin{bmatrix} \mathbb{C}^4 \\ \mathbb{C}^2 \end{bmatrix}
 \to  \begin{bmatrix} \mathbb{C}^4   \\ \mathbb{C}^2 \end{bmatrix} , \\
B_w & =  \begin{bmatrix}  B_4 & 0 & 0_{ 4 \times 3 }  \\
                        0 &  B_2 & 0_{ 2 \times 3 }  \end{bmatrix}:
\begin{bmatrix}  \mathbb{C} \\  \mathbb{C}  \\  \mathbb{C}^3 \end{bmatrix}  \to \begin{bmatrix} \mathbb{C}^4 \\ \mathbb{C}^2 \end{bmatrix} , \\
 C_w &=
\begin{bmatrix} C_4 & 0 \\ 0 & C_2 \\ 0_{3 \times 4} &0_{3 \times 2} \end{bmatrix}:
 \begin{bmatrix}  \mathbb{C}^4 \\  \mathbb{C}^2 \\  \end{bmatrix}\rightarrow \begin{bmatrix}  \mathbb{C} \\  \mathbb{C}  \\  \mathbb{C}^3  \end{bmatrix}
 = \BC^5 , \\
 D_w &= \begin{bmatrix}   D_4 & 0 & 0 \\  0 & D_2 & 0 \\ 0 & 0 & I_3 \end{bmatrix}:
\begin{bmatrix}  \mathbb{C} \\  \mathbb{C}  \\  \mathbb{C}^3  \end{bmatrix}\to   \begin{bmatrix}  \mathbb{C} \\  \mathbb{C}  \\  \mathbb{C}^3  \end{bmatrix} .
\end{align*}

It is noted that $B_vB_w^* = 0$. Therefore the
unique solution to the Lyapunov equation
$A_v \Omega +\Omega A_w^* + B_v B_w^*=0$
is $\Omega =0$. Hence
\begin{align*}
C_\circ &=  D_v B_w^* + C_v \Omega = D_v B_w^*   =
  \sqrt{2} \begin{bmatrix}
     -1  &  1 &  -1  &  1   &      0   &      0 \\
         0    &     0   &      0    &     0  & -1  &  1 \\
         0    &     0   &      0   &      0  &       0   &      0  \\
         0    &     0    &     0    &     0   &      0    &     0  \\
         0    &     0    &     0    &     0   &      0    &     0  \\
   \end{bmatrix}.
\end{align*}
The unique solution to the Lyapunov equation
$A_w  Q +Q  A_w^* + C_\circ^* C_\circ =0$ is given by $Q = I$.
Now observe that
\[
\zeta(-A_w) = \begin{bmatrix}
                 0  &   1  &    0   &   0  &    0 &     0 \\
                 0  &   0  &    1   &   0  &    0 &     0 \\
                 0  &   0  &    0   &   1  &    0 &     0 \\
                 0  &   0  &    0   &   0  &    0 &     0 \\
                 0  &   0  &    0   &   0  &    0 &     1 \\
                 0  &   0  &    0   &   0  &    0 &     0 \\
          \end{bmatrix}.
          \]
It is noted that $\zeta(-A_w) =S_4 \oplus S_2$ on $\mathbb{C}^4 \oplus \mathbb{C}^2$ where
$S_k$ on $\mathbb{C}^k$ is the upward shift, that is, all the entries of $S_k$
immediately above the main diagonal are $1$ and all the other entries are zero,
that is $ S_k = J_k(0) $. (See also Corollary \ref{CtoD}.)
Because $Q=I$, we readily see that
\begin{align*}
& \dim (\kr (I-Q)) = 6\\
&\dim(\kr(I- \zeta(-A_w)  Q \zeta(-A_w^*) ) )= 4 \\
& \dim\left(\kr \left(I - \zeta(-A_w)^2 Q \zeta(-A_w^*)^{2} \right) \right)= 2\\
&\dim\left(\kr \left(I - \zeta(-A_w)^3 Q \zeta(-A_w^*)^{3} \right)\right) = 1 \\
& \dim\left(\kr \left(I - \zeta(-A_w)^k Q \zeta(-A_w^*)^k \right) \right)= 0 \mbox{\ \ for  } k \geq 4.
\end{align*}
In other words,
\[
\mu_1 = 6-4=2, \quad \mu_2 = 4-2=2,\quad \mu_3 = 2-1=1 \mbox{ and } \mu_4=1-0=0.
\]
Using this we have
\[
\kappa_2 = \#\{k:\mu_k \geq 1\} = 4
\quad \mbox{and}\quad \kappa_1 = \#\{k:\mu_k \geq 2\} = 2.
\]
Therefore the negative Wiener-Hopf indices for $R$ are $\{-4,-2\}$.


\subsection*{Acknowledgements}

The work of the second author is supported in part by the National Research Foundation of South Africa (NRF), Grant Number 145688.


\end{document}